\documentclass[12pt,a4paper,notitlepage,openbib]{article}
[2000/05/19 v1.4b Standard LaTeX document class]

\usepackage{amsmath}
\usepackage{amsfonts}
\usepackage{amssymb}
\usepackage{amsthm}
\usepackage{amsrefs}

\usepackage[cp1251]{inputenc}
\usepackage[english]{babel} 

\usepackage[dvips]{graphicx}
\usepackage[dvips]{color}
\usepackage[english, xnumber]{automata}

\author{L. Bartholdi\and I. I. Reznykov\and V. I. Sushchansky}

\title{The smallest Mealy automaton of intermediate growth}

\date{2004/02/10}

\begin{document}

\maketitle

\tableofcontents

\begin{abstract}
  In this paper we study the automaton $\Intermediate{2}$, the
  smallest Mealy automaton of intermediate growth, first considered
  in~\cite{reznykov-s:interm2x2}. We describe the automatic
  transformation monoid defined by $\Intermediate{2}$, give a formula
  for the generating series for the (ball volume) growth function of
  $\Intermediate{2}$, and give sharp asymptotics for the growth
  function of $\Intermediate{2}$, namely
  \[
    \Growth{\Intermediate{2}} \n \sim 2 ^{5/2} 3 ^{3/4} \pi ^{-2} n^{1/4}
    \exp{\Argument{\pi \sqrt{{n}/{6}}}},
  \]
  with the ratios of left- to right-hand side tending to $1$ as
  $n \to \infty$.
\end{abstract}


\section{Introduction}

\par The growth of a Mealy automaton is defined as
the growth of the number of pairwise inequivalent internal states of
iterates of that automaton. This notion of growth was introduced by
R.I.~Grigorchuk in~\cite{MR89f:20065}.  The growth function of an
arbitrary Mealy automaton coincides with the spherical growth function
of the automatic transformation semigroup it defines, and actually the
growth of automata are calculated by investigating the growth of the
corresponding automatic transformation semigroups.

The automatic transformation groups defined by invertible $2$-state
Mealy automata over the $2$-symbol alphabet were described
in~\cite{grigorchuk-n-s:automata}. The automatic transformation
semigroups defined by all $2$-state Mealy automata over the $2$-symbol
alphabet were investigated in~\cites{reznykov:phd} and in the
papers~\cites{reznykov-s:growths,reznykov-s:fibonacci,reznykov-s:interm2x2}.

\par Among these semigroups there are twelve finite semigroups, seven
semigroups of polynomial growth, one semigroup of intermediate growth,
and eight semigroups of exponential growth, including the free
semigroup. There are four pairwise similar (in the sense of
Definition~\ref{def:similar_automata}) $2$-state Mealy automata over
the $2$-symbol alphabet of intermediate growth order, and these
automata define isomorphic automatic transformation semigroups. One of
these automata was considered in~\cite{reznykov:phd} and~\cite{reznykov-s:interm2x2}.  There, an automatic transformation
semigroup of intermediate growth was constructed, with an exact
formula for the growth function, expressed as an infinite sum. Its
growth order was estimated between $\GrowthOrder{ e ^{\sqrt[4]{n} \,
  }}$ and $\GrowthOrder{ e ^{\sqrt{n} \, }}$.

\par In this paper we consider the automaton of intermediate growth
$\Intermediate{2}$ and the semigroup of automatic transformations
$\Semigroup{\Intermediate{2}}$ that it defines. In
Theorem~\ref{th:semigroup} we describe the semigroup
$\Semigroup{\Intermediate{2}}$ and its quotient semigroups, in
Theorem~\ref{th:generating_functions} we exhibit the growth series of
the automaton and the semigroup, and in Theorem~\ref{th:estimates} we
derive sharp asymptotics for the growth functions. The first part of
Theorem~\ref{th:semigroup} was proved in~\cite{reznykov:phd} and~\cite{reznykov-s:interm2x2}, but we give here a shorter proof, and
a new proof of the minimality of the system of defining relations.
Moreover, the other results are new.

\par There are various motivations for the precise study of growth
functions of semigroups generated by automata. The first, and in some
sense only, known examples of groups of intermediate growth come from
automata~\cite{grigorchuk:growth}, and these groups' structure can at
least partly be understood through their growth.  Also, the natural
algebraic object associated to a Mealy automaton is a semigroup, which
is a group only under an additional assumption.  Furthermore, it seems
beyond reach to obtain as sharp results as those of this paper for
even the simplest known groups of intermediate growth.

\par Finally, a word should be added as to what is meant by deriving
an ``exact formula'' for the growth of a semigroup, that is not
tautological. The formulae we obtain in this paper have the merits of
being easily and quickly computable, and of being expressible
algebraically in terms of the partition function. This is certainly
the most that can be hoped from a transcendental generating series.

\section{Main results}

\begin{figure}[t]
  \centering
  \includegraphics*[trim = 1 1 1 1, height = 7cm, width = 9cm, keepaspectratio = true]{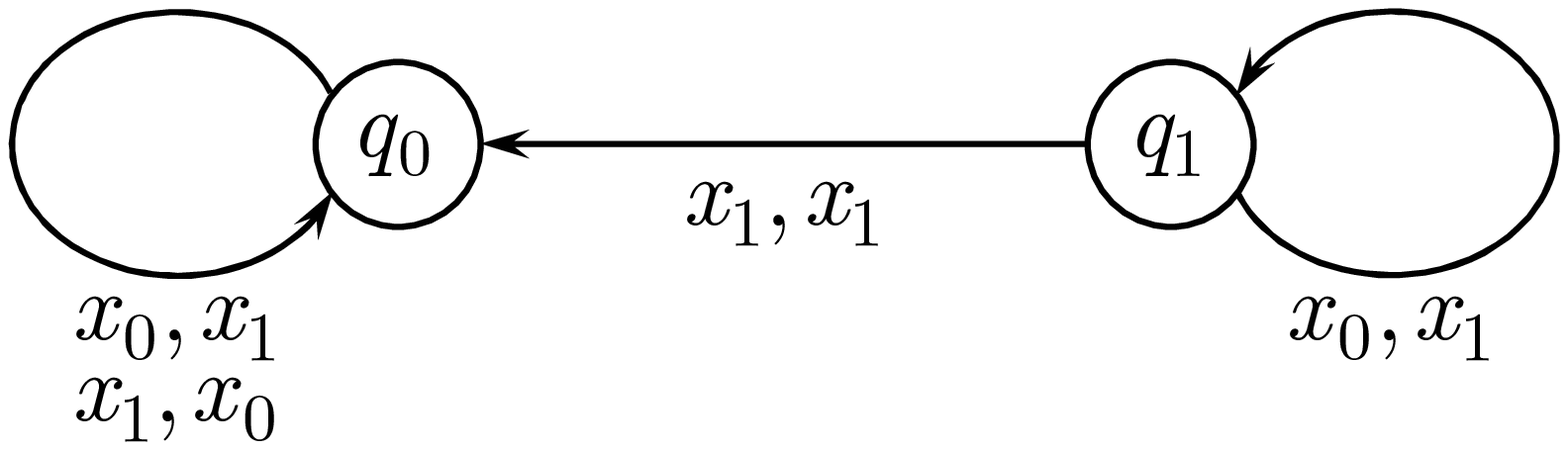}
  \caption{The automaton $\Intermediate{2}$}
  \label{fig:automaton_intermediate}
\end{figure}

\par Let $\Intermediate{2}$ be the $2$-state Mealy automaton over the
$2$-symbol alphabet whose Moore diagram is shown on
Figure~\ref{fig:automaton_intermediate}. Let us denote the semigroup
defined by $\Intermediate{2}$ by the symbol
$\Semigroup{\Intermediate{2}}$, and the growth functions of
$\Intermediate{2}$ and $\Semigroup{\Intermediate{2}}$ by the symbols
$\Growth{\Intermediate{2}}$ and $\GrowthSemigroup{\Intermediate{2}}$,
respectively. Let us denote for each $\PositiveInteger$ the quotient
semigroup given by the representation of $\Intermediate{2}$ as maps
from $\{x_0,x_1\}^n$ to itself by the symbol $\Quotient$. The
following theorem holds:

\begin{theorem} \label{th:semigroup}
\begin{enumerate}
\item The semigroup $\Semigroup{\Intermediate{2}}$ is a monoid, and
  has the following
  presentation~\cites{reznykov-s:interm2x2,reznykov:phd}:
\begin{equation} \label{eq:semigroup}
    \Semigroup{\Intermediate{2}} = \ATSemigroup{f_0 ^2  = 1;\,f_1 \left( {f_0
    f_1} \right) ^p \left( {f_1 f_0} \right) ^p f_1 ^2  = f_1 \left( {f_0 f_1 }
    \right) ^p \left( {f_1 f_0 } \right) ^p, \, p \ge 0}.
\end{equation}
The monoid $\Semigroup{\Intermediate{2}}$ is infinitely presented, and
the word problem is solvable in polynomial time.
  \item The semigroup $\Quotient$, $\PositiveInteger$, has the presentation \[
    \Quotient = \ATSemigroup{
        \begin{array}{*{20}l}
            {f_0 ^2  = 1}; \\
            {f_1 \left( {f_0 f_1} \right) ^p \left( {f_1 f_0} \right) ^p f_1 ^2
            = f_1 \left( {f_0 f_1 } \right) ^p \left( {f_1 f_0 } \right) ^p, \,
            0 \le p \le n - 2}; \\
            {f_1 \left( {f_0 f_1} \right) ^{n - 1} f_1 = f_1 \left( {f_0 f_1}
            \right) ^{n - 1} f_0 = f_1 \left( {f_0 f_1} \right) ^{n - 1}}
        \end{array}}
    \]
\end{enumerate}
\end{theorem}

\par The following corollary follows (for relevant definitions
see Section~\ref{subsect:Hausdorff_dimension}):
\begin{corollary} \label{cor:Hausdorff_dimension}
  The semigroup $\Semigroup{\Intermediate{2}}$ has Hausdorff dimension $0$.
\end{corollary}

\begin{theorem} \label{th:generating_functions}
\begin{enumerate}
\item The word growth series $\Delta _{\Semigroup{\Intermediate{2}}} \Argument{X} = \sum \limits_{n \ge 0}
  {\WordGrowthSemigroup{\Intermediate{2}} \n X ^n}$ of
  $\Semigroup{\Intermediate{2}}$ admits the description
\[
    \Delta _{\Semigroup{\Intermediate{2}}} \Argument{X} = \Argument{1 + X}
    \left( 1 + \frac{X}{1 - X} \prod \limits_{n \ge 0} \Argument{1 + X ^{2n +
    1}} \right).
\]

\item The growth series $\Gamma _{\Intermediate{2}} \Argument{X} =
  \sum \limits_{n \ge 0} {\Growth{\Intermediate{2}} \n X ^n}$ of
  $\Intermediate{2}$ admits the description
\[
    \Gamma _{\Intermediate{2}} \Argument{X} = \frac{1}{1 - X} \left( 1 +
    \frac{X}{1 - X} \prod \limits_{n \ge 0} \Argument{1 + X ^{2n + 1}} \right).
\]

\item The growth series $\Gamma _{\Semigroup{\Intermediate{2}}}
  \Argument{X} = \sum \limits_{n \ge 0}
  {\GrowthSemigroup{\Intermediate{2}} \n X ^n}$ of
  $\Semigroup{\Intermediate{2}}$ admits the description
\[
    \Gamma _{\Semigroup{\Intermediate{2}}} \Argument{X} = \frac{1 + X}{1 - X}
    \left( 1 + \frac{X}{1 - X} \prod \limits_{n \ge 0} \Argument{1 + X ^{2n +
    1}} \right).
\]

\end{enumerate}
\end{theorem}

\par Let us denote the number of all partitions of a positive integer
$n$ to $k$ odd parts by the symbols $\OddPartitions{n}$.

\begin{theorem} \label{th:estimates}
The growth functions have the following sharp estimates:
\begin{align*}
    \WordGrowthSemigroup{\Intermediate{2}} \n & \sim \frac{4 \sqrt{6}}{\pi}
    \sqrt{n} \cdot \OddPartitions{n} \sim \frac{2 ^2 3 ^{1/4}}{\pi} n^{-1/4}
    \exp{\Argument{\pi \sqrt{\frac{n}{6}}}};\\
    \Growth{\Intermediate{2}} \n & \sim \frac{24}{\pi ^2} n \cdot
    \OddPartitions{n} \sim \frac{2 ^{5/2} 3 ^{3/4}}{\pi ^2} n^{1/4}
    \exp{\Argument{\pi \sqrt{\frac{n}{6}}}};\\
    \GrowthSemigroup{\Intermediate{2}} \n & \sim \frac{48}{\pi ^2} n \cdot
    \OddPartitions{n} \sim \frac{2 ^{7/2} 3 ^{3/4}}{\pi ^2} n^{1/4}
    \exp{\Argument{\pi \sqrt{\frac{n}{6}}}}.\\
\end{align*}
\end{theorem}

\begin{corollary} \label{cor:growth_orders}
The growth orders of the growth
  functions of $\Intermediate{2}$ and $\Semigroup{\Intermediate{2}}$
  are equal, and
\[
    \GrowthOrder{\Growth{\Intermediate{2}}} =
    \GrowthOrder{\GrowthSemigroup{\Intermediate{2}}} = \GrowthOrder{\exp
    \Argument {\sqrt {n}}}.
\]
\end{corollary}

\section{Preliminaries}

By $\Natural$ we mean the set of non-negative integers
$\Natural=\{0,1,2,\dots\}$.

\subsection{Growth functions}

\par Let us consider the set of positive non-decreasing functions of
a natural argument $\Growth{} : \Natural \to \Natural$; in the sequel
such functions will be called \emph{growth functions}.

\begin{definition} \label{def:growth_order}
  For $i = 1, 2$ let $\Growth{i} : \Natural \to \Natural$ be growth
  functions. The function $\Growth{1}$ has \textit{no greater growth
    order} (notation $\Growth{1} \preceq \Growth{2}$) than the
  function $\Growth{2}$, if there exist numbers $C_1, C_2, N_0 \in
  \Natural$ such that
  \[\Growth{1} \n \le C_1 \Growth{2} \Argument{C_2 n}\]
  for any $n \ge N_0$.
\end{definition}

\begin{definition} \label{def:growth_equivalence}
  The growth functions $\Growth{1}$ and $\Growth{2}$ are equivalent or
  have \textit{the same growth order} (notation $\Growth{1} \sim
  \Growth{2}$), if the following inequalities hold: \[\Growth{1}
  \preceq \Growth{2} \quad \text{and} \quad \Growth{2} \preceq
  \Growth{1}.\]
\end{definition}

\par The equivalence class of the function $\Growth{}$ is called
its \textit{growth order} and is denoted by the symbol
$\GrowthOrder{\Growth{}}$. The relation $\preceq$ induces an order relation,
written $<$, on equivalence classes. The growth order $\GrowthOrder{\Growth{}}$
is called \textit{intermediate} if $\GrowthOrder{n ^d} <
\GrowthOrder{\Growth{}} < \GrowthOrder{e ^n}$ for any $d > 0$.

\subsection{Mealy automata}

\par For $m\ge 2$ let $\Alphabet{m}$ be the $m$-symbol alphabet,
$\Alphabet{m} = \left\{ x_0, x_1, \dots, x_{m - 1} \right\}$. Let us
denote the set of all finite words over $\Alphabet{m}$, including the
empty word $\varepsilon$, by the symbol $\FiniteSet$, and denote the
set of all infinite (to right) words by the symbol $\InfiniteSet$.

\par Let $\Automaton = \MAutomaton{n}{m}{}$ be a \textit{non-initial
  Mealy automaton}~\cite{MR17:436b} with finite set of states
$\States{n} = \left\{ q_0, q_1, \dots, q_{n - 1} \right\}$; input and
output alphabets are the same and are equal to $\Alphabet{m}$, and
$\pi : \Alphabet{m} \times \States{n} \to \States{n}$ and $ \lambda :
\Alphabet{m} \times \States{n} \to \Alphabet{m}$ are its transition
and output functions, respectively.


\par The function $\lambda$ can be extended in a natural way to a
mapping $\lambda: \FiniteSet \times \States{n} \to \FiniteSet$ or to
a mapping $\lambda: \InfiniteSet \times \States{n} \to \InfiniteSet$
(see for example,~\cite{glushkov:ata}).

\begin{definition} \label{def:automatic_transformation}
  For any state $\q \in \States{n}$ the transformation $f _{\q,
    {\Automaton}} : \FiniteSet \to \FiniteSet$, respectively $f _{\q,
    {\Automaton}} : \InfiniteSet \to \InfiniteSet$, defined by
\[f _{\q, {\Automaton}} \Argument{u} = \lambda \Argument{u, \q},\] where $u
\in \FiniteSet$, respectively $u \in \InfiniteSet$, is called the
\textit{automatic transformation} defined by $\Automaton$ at the state
$\q$.
\end{definition}

\par We write  a function $f:\Alphabet{m}\to\Alphabet{m}$ as $\left(
{f\Argument{x_0}}, {f\Argument{x_1}}, \ldots, {f\Argument{x_{m - 1}}} \right)$.
Let us consider the transformation $\sigma _ \q$ over the alphabet
$\Alphabet{m}$, $\q \in \States{n}$, defined by the output function $\lambda$:
\[
   \sigma _{\q} = \left(
   {\lambda \Argument{x_0, \q}}, \, {\lambda \Argument{x_1, \q}}, \,
   \ldots, \, {\lambda \Argument{x _{m - 1} ,\q}}
\right). \]

\par
Interpreting an automatic transformation as an isomorphism of the
rooted $m$-regular tree (see for
example~\cite{grigorchuk-n-s:automata}), we have the
following interpretation. Let $\q$ be an arbitrary state. The image of
the word $u = u_0 u_1 u_2 \cdots \in \InfiniteSet$ under the action of
the automatic transformation $f _{\q, \Automaton}$ can be written in
the following way:
\[
f _{\q, \Automaton} \Argument{u_0 u_1 u_2 \dots} = \lambda
\Argument{u_0, \q} \cdot f _{\pi \Argument{u_0, \q}, \Automaton}
\Argument{u_1 u_2 \dots } = \sigma _{\q} \Argument{u_0} \cdot f _{\pi
  \Argument{u_0, \q}, \Automaton} \Argument {u_1 u_2 \dots}.
\]
This means that $f _{\q, \Automaton}$ acts on the first symbol of the
word $u$ by the transformation $\sigma _{\q}$ over the alphabet
$\Alphabet{m}$, and acts on the remainder of the word without its
first symbol by the transformation $f _{\pi \Argument{u_0, \q},
  \Automaton}$. Therefore the transformations defined by the automaton
$\Automaton$ can be written in \textit{unrolled form}:
\[ f _{q_i} = \left( {f _{\pi \Argument{x_0 ,q_i}}, f _{\pi
      \Argument{x_1, q_i}}, \dots, f _{\pi \Argument{x_{m - 1}, q_i}}
  } \right) \sigma _{q_i}, \] where $i = 0, 1, \dots, {n - 1}$.

\par Let us illustrate this notion. Let $\Intermediate{2}$ be the
automaton, shown on Figure~\ref{fig:automaton_intermediate}, and let us
construct the unrolled forms of its automatic transformations. As $\pi
\Pair{x_0}{q_0} = \pi \Pair{x_1}{q_0} = q_0$ and $\sigma _{q_0} = \TX{0}{1}$,
then the unrolled form of $f_{q_0}$ is written as
\[
    f_{q_0} = \Argument{f_{q_0}, f_{q_0}} \TX{0}{1}.
\]
Similarly, there are $\pi \Pair{x_0}{q_1} = q_1$, $\pi \Pair{x_1}{q_1} = q_1$
and $\sigma _{q_1} = \TX{1}{1}$. Hence the unrolled form of $f_{q_1}$ is
defined by
\[
    f_{q_1} = \Argument{f_{q_1}, f_{q_0}} \TX{1}{1}.
\]
Let $u = x_0 x_0 x_1 x_0 x_0 x_1 \ldots = \left( x_0 x_0 x_1 \right)
^\ast$ be an infinite word, and let us consider the action of
$f_{q_0}$ and $f_{q_1}$ on it.  We have
\begin{align*}
    f_{q_0} \Argument{u} & = \sigma _{q_0} \Argument{x_0} \cdot f_{q_0}
    \Argument{x_0 x_1 x_0 x_0 x_1 \ldots} = x_1 \cdot \sigma _{q_0}
    \Argument{x_0} \cdot f_{q_0} \Argument{x_1 x_0 x_0 x_1 \ldots} =\\
    & = x_1 x_1 \cdot \sigma _{q_0} \Argument{x_1} \cdot f_{q_0} \Argument{x_0
    x_0 x_1 \ldots} = x_1 x_1 x_0 \cdot f_{q_0} \Argument{u} = \ldots = \\
    & = x_1 x_1 x_0 x_1 x_1 x_0 \ldots = \left( x_1 x_1 x_0 \right) ^\ast,
\end{align*}
and
\begin{align*}
    f_{q_1} \Argument{u} & = \sigma _{q_1} \Argument{x_0} \cdot f_{q_1}
    \Argument{x_0 x_1 x_0 x_0 x_1 \ldots} = x_1 \cdot \sigma _{q_1}
    \Argument{x_0} \cdot f_{q_1} \Argument{x_1 x_0 x_0 x_1 \ldots} =\\
    & = x_1 x_1 \cdot \sigma _{q_1} \Argument{x_1} \cdot f_{q_0} \Argument{x_0
    x_0 x_1 \ldots} = x_1 x_1 x_1 \cdot f_{q_0} \Argument{u} = \\
    & = x_1 x_1 x_1 \cdot \left( x_1 x_1 x_0 \right) ^\ast.
\end{align*}

\par
The Mealy automaton $\Automaton = \MAutomaton{n}{m}{}$ defines the set $F
_{\Automaton} = \left\{ f _{q_0}, f _{q_1 }, \dots, f _{q_{n - 1}} \right\}$ of
automatic transformations over $\FiniteSet$. The Mealy automaton $\Automaton$
is called \textit{invertible} if all transformations from the set $F
_{\Automaton}$ are bijections. It is easy to show (see for
example~\cite{grigorchuk-n-s:automata}) that $\Automaton$ is
invertible if and only if the transformation $\sigma _{\q}$ is a permutation of
$\Alphabet{m}$ for each state $\q \in \States{n}$.

\begin{definition}[\cite{glushkov:ata}]\label{def:automata_isomorphic}
  The Mealy automata $\Automaton[i] = \MAutomaton{n}{m}{i}$ for $i = 1,
  2$ are called \textit{isomorphic} if there exist permutations $\xi,
  \psi \in Sym \Argument{\Alphabet{m}}$ and $\theta \in Sym
  \Argument{\States{n}}$ such that
\begin{align*}
\theta \pi _1 \Argument{\x, \q} & = \pi _2 \Argument{\xi \x, \theta \q}, & \psi
\lambda _1 \Argument{\x, \q} & = \lambda _2 \Argument{\xi \x,\theta \q}
\end{align*}
for all $\x \in \Alphabet{m}$ and $\q \in \States{n}$.
\end{definition}

\begin{definition}[\cite{glushkov:ata}]
  The Mealy automata $\Automaton[i] = \MAutomaton{n_i}{m}{i}$ for $i = 1,
  2$, are called \textit{equivalent} if $F _{\Automaton[1]} = F
  _{\Automaton[2]} $.
\end{definition}

\begin{proposition}[\cite{glushkov:ata}]
  Each class of equivalent Mealy automata over the alphabet
  $\Alphabet{m}$ contains, up to isomorphism, a unique automaton that
  is minimal with respect to the number of states (such an automaton
  is called \textit{reduced}).
\end{proposition}
The minimal automaton can be found using the standard algorithm of
minimization.

\begin{definition}[\cite{gecseg:products}] \label{def:automata_product}
  For $i=1,2$ let $\Automaton[i] = \MAutomaton{n_i}{m}{i}$ be
  arbitrary Mealy automata. The automaton $\Automaton = \left(
    \Alphabet{m}, \States{n _1} \times \States{n _2}, \pi, \lambda
  \right)$ such that its transition and output functions are defined
  in the following way:
\begin{align*}
  \pi \Argument{\x, \Pair{\q_1}{\q_2}} & = \Pair{\pi _1
    \Argument{\lambda _2
      \Argument{\x, \q_2}, \q_1}}{\pi _2 \Argument{\x, \q_2}},\\
  \lambda \Argument{\x, \Pair{\q_1}{\q_2}} & = \lambda _1
  \Argument{\lambda _2 \Argument{\x, \q_2}, \q_1 },
\end{align*}
where $\x \in \Alphabet{m}$ and $\Pair{\q_1}{\q_2} \in \States{n _1}
\times \States{n _2}$, is called the \emph{product} of the automata
$\Automaton[1]$ and $\Automaton[2]$.
\end{definition}

\begin{proposition}[\cite{gecseg:products}] \label{prop:automaton_multiple}
  For any states $\q_1 \in \States{n _1}$ and $\q_2 \in \States{n _2}$
  and an arbitrary word $u \in \FiniteSet$ the following equality
  holds:
\[ f _{\Pair{\q_1}{\q_2}, \Automaton} \Argument{u} = f _{\q_1, \Automaton[1]}
\Argument{f _{\q_2, \Automaton[2]} \Argument{u}}. \]
\end{proposition}

\par It follows from Proposition~\ref{prop:automaton_multiple} that
for the transformations $f _{\q_1, \Automaton[1]}$ and $f _{\q_2,
  \Automaton[2]}$, with $\q_1 \in \States{n _1}$ and $\q_2 \in
\States{n _2}$, the unrolled form of the product $f _{\Argument{\q_1,
    \q_2}, {\Automaton[1] \times \Automaton[2]}}$ is defined by:
\[ f _{\Argument{\q_1, \q_2}, {\Automaton[1] \times \Automaton[2]}} = f
_{\q_1, \Automaton[1]} {f _{\q_2, \Automaton[2]}} = \left( g_0, g_1, \dots,
g_{m - 1} \right) \sigma _{\q_1, \Automaton[1]} \sigma _{\q_2,
\Automaton[2]},\] where $g _i = f _{\pi_1 \Argument{\sigma _{\q_2
,\Automaton[2]} \Argument{x_i}, \q_1}, \Automaton[1]} f _{\pi_2 \Argument{x_i,
\q_2}, \Automaton[2]}$ for $i = 0, 1, \dots, {m - 1}$.

\par The power $\Automaton ^n$ is defined for any automaton
$\Automaton$ and any positive integer $n$. Let us denote $\Automaton
^{\n}$ the minimal Mealy automaton equivalent to $\Automaton ^n$. It
follows from Definition~\ref{def:automata_product} that $\left|
  {\States{\Automaton ^{\n}}} \right| \le \left| {\States{\Automaton}}
\right| ^n$.

\begin{definition}[\cite{MR89f:20065}] \label{def:growth_automaton}
  The function $\GrowthAutomaton{}$ of a natural argument $n\ge1$,
  defined by
  \[\GrowthAutomaton{} \n = \left| {\States{\Automaton^{\n}}} \right|\]
  is called the \textit{growth function} of the Mealy automaton
  $\Automaton$.
\end{definition}

\begin{definition}[\cite{reznykov:phd}] \label{def:similar_automata}
  The Mealy automata $\Automaton[i] = \MAutomaton{n}{m}{i}$, for $i =
  1, 2$, are called \textit{similar} if they are isomorphic in the
  sense of Definition~\ref{def:automata_isomorphic}, for permutations
  $\xi, \psi \in Sym \Argument{\Alphabet{m}}$ satisfying furthermore
  $\psi = \xi$.
\end{definition}

\subsection{Semigroups}

\par Let $\Semigroup{}$ be a semigroup with the finite set of
generators $G = \left\{ s_0, s_1, \dots, s_{k - 1} \right\}$. Let us
denote the free semigroup with the set $G$ of generators by the symbol
$G^{+}$. It is easy to see (for example,
in~\cite{lallement:semigroups}) that if the semigroup $\Semigroup{}$
does not contain the identity, then $\Semigroup{}$ is a homomorphic
image of the free semigroup $G^{+}$. Similarly, the monoid
$\Semigroup{} = \mathop{sg} \Argument{G}$ is a homomorphic image of
the free monoid $G^{*}$.

\par The elements of the free semigroup $G^{+}$ are called
\textit{semigroup words}. In the sequel, we identify them with
corresponding elements of $\Semigroup{}$. The semigroup words $\s_1$
and $\s_2$ are called \textit{equivalent relative to the system $G$ of
  generators in the semigroup $\Semigroup{}$}, if in $\Semigroup{}$
the equality $\s_1 = \s_2$ holds~\cite{lallement:semigroups}.

\begin{definition} \label{def:length_element}
  Let $\s$ be an arbitrary element of $\Semigroup{}$. The length
  $\Length{\s}$ of $\s$ is the minimal possible number $\ell>0$ of
  generators in a factorization \[\s = {s_{i_1}} {s_{i_2}} {s_{i_3}}
  \dots {s_{i_\ell}},\] where $s_{i_j} \in G$ for all $1 \le j
  \le\ell$.
\end{definition}

\par Obviously for any $\s \in \Semigroup{}$ the length $\Length{\s}$
is greater than $0$; but let us assume $\Length{1} = 0$, if
$\Semigroup{}$ is a monoid.

\par Let us order the generators of $\Semigroup{}$ according to their
index; and introduce a linear order on the set of elements of $G^{+}$:
semigroup words are ranked by length, and then words of the same
length are arranged lexicographically. The \emph{representative} of a
class in the equivalence relation introduced above is the minimal
semigroup word in the sense of this order.

\begin{definition}
  Let $\s \in \Semigroup{}$ be an arbitrary element. The \emph{normal
    form} of this element is the representative of the equivalence
  class of semigroup words mapped to the element $\s$.
\end{definition}

\begin{definition} \label{def:growth_semigroup}
  The function $\GrowthSemigroup{}$ of a natural argument
  $\PositiveInteger$ defined by
\[ \GrowthSemigroup{} \n =
\left| \left\{ {\begin{array}{*{20}c}
  {s \in \Semigroup{}} & \vline & {\Length{s} \le n}  \\
\end{array}} \right\} \right|\]
is called the \textit{growth function of $\Semigroup{}$ relative to
  the system $G$ of generators}.
\end{definition}

\begin{definition} \label{def:spherical_growth_semigroup}
  The function $\SphericalGrowthSemigroup{}$ of a natural argument
  $\PositiveInteger$ defined by
\[ \SphericalGrowthSemigroup{} \n =
\left| \left\{ {\begin{array}{*{20}c}
  {s \in \Semigroup{}} & \vline & {s = s_{i_1} s_{i_2} \dots s_{i_n}, \,
  s_{i_j} \in G, \, 1 \le j \le n}  \\
\end{array}} \right\} \right|
\] is called the \textit{spherical growth function of $\Semigroup{}$ relative
to the system $G$ of generators}.
\end{definition}

\begin{definition} \label{def:word_growth_semigroup}
  The function $\WordGrowthSemigroup{}$ of a natural argument
  $\PositiveInteger$ defined by
  \[\WordGrowthSemigroup{} \n =
  \left| \left\{ {\begin{array}{*{20}c}
          {s \in \Semigroup{}} & \vline & {\Length{s} = n}  \\
        \end{array}} \right\} \right|
  \] is called the \textit{word growth function of $\Semigroup{}$
    relative to the system $G$ of generators}.
\end{definition}

If we denote by $\pi: G ^{+} \to \Semigroup{}$ the natural epimorphism from the
free semigroup $G^ {+}$ to $\Semigroup{}$, these functions can be expressed as
follows:
\begin{align*}
  \GrowthSemigroup{} \n &= \Card{\bigcup_{i = 0} ^n \pi \Argument{G^i}},\\
  \SphericalGrowthSemigroup{} \n &= \Card{\pi \Argument{G ^n}},\\
  \WordGrowthSemigroup{} \n &= \Card{\pi \Argument{G^n} \setminus \bigcup _{i =
  0} ^{n - 1} \pi \Argument{G ^i}}.
\end{align*}

\par The following proposition is well-known, and is proved in many
papers (see for
example~\cites{grigorchuk-n-s:automata,nathanson:amemm}):

\begin{proposition} \label{prop:equivalent_growth_semigroup}
  Let $\Semigroup{}$ be an arbitrary finitely generated semigroup, and
  let $G_1$ and $G_2$ be systems of generators of $\Semigroup{}$.  Let
  us denote the growth function of $\Semigroup{}$ relative to the set
  $G_i$ of generators by the symbol $\GrowthSemigroup{i}$, for $i = 1,
  2$. Then $\GrowthOrder{\GrowthSemigroup{1}} =
  \GrowthOrder{\GrowthSemigroup{2}}$.
\end{proposition}

\par From
Definitions~\ref{def:growth_semigroup},~\ref{def:spherical_growth_semigroup}
and~\ref{def:word_growth_semigroup}, the following inequalities hold
for any $\PositiveInteger$:
\begin{equation} \label{eq:growths_semigroup}
\WordGrowthSemigroup{} \n \le \SphericalGrowthSemigroup{} \n \le
\GrowthSemigroup{} \n = \sum \limits _{i = 0} ^n {\WordGrowthSemigroup{}
\Argument{i}}.
\end{equation}

\begin{proposition} \label{prop:monoid_growth_order}
  Let $\Semigroup{}$ be an arbitrary finitely generated monoid. Then
\[
    \GrowthOrder{\WordGrowthSemigroup{}} \le
    \GrowthOrder{\SphericalGrowthSemigroup{}} =
    \GrowthOrder{\GrowthSemigroup{}}.
\]
\end{proposition}

\par Let $\Semigroup{}$ be a semigroup without identity. Then the
growth function and the spherical growth function may have different
growth orders.  For example, let $\Semigroup{} = \Natural$ be the
additive semigroup $\Semigroup{} = \mathop{sg} \Argument{1}$. Then
$\GrowthSemigroup{} \n = n$, $\SphericalGrowthSemigroup{} \n = 1$, and
these functions have different growth orders, $\GrowthOrder{1} <
\GrowthOrder{n}$.

\par There are many results concerning the growth of groups. For
references see the survey~\cite{grigorchuk-n-s:automata}, or the
book~\cite{harpe:ggt}.

\subsection{Growth series}

\par It is often convenient to encode the growth function of a semigroup in a
generating series:
\begin{definition}
  Let $\Semigroup{}$ be a semigroup generated by a finite set $G$. The
  \emph{growth series} of $\Semigroup{}$ is the formal power series \[
  \Gamma _{\Semigroup{}} \n[X] = \sum \limits _{n \ge 0}
  \GrowthSemigroup{} \n X ^n.
  \]
\end{definition}
The power series $\Delta _{\Semigroup{}} \n[X] = \sum \limits_{n \ge 0}
\WordGrowthSemigroup{} \n X ^n$ can also be introduced; we then have $\Delta
_{\Semigroup{}} \n[X] = \Argument{1 - X} \Gamma _{\Semigroup{}} \n[X]$. The
series $\Delta _{\Semigroup{}}$ is called the \textit{word growth series} of
the semigroup $\Semigroup{}$.

\par The growth series of a Mealy automaton is introduced similarly:
\begin{definition}
  Let $\Automaton$ be an arbitrary Mealy automaton. The \emph{growth
    series} of $\Automaton$ is the formal power series \[ \Gamma
  _{\Automaton} \n[X] = \sum \limits _{n \ge 0} \GrowthAutomaton{} \n
  X ^n.
  \]
\end{definition}

\par The radius of convergence, and behaviour of $\Gamma _{\Semigroup{}}$ near
its singularities, encode the asymptotics of $\GrowthSemigroup{}$. The
semigroup $\Semigroup{}$ has subexponential growth if and only if $\Gamma
_{\Semigroup{}}$ converges in the open unit disk.

\par Sharper results of this flavour are often called \emph{tauberian} and
\emph{abelian} theorems.  We quote two such results~\cite{nathanson:density}:
\begin{theorem}
  If $\Gamma _{\Semigroup{}}$ converges in the open unit disk, and
  $\log \GrowthSemigroup{} \n \sim 2\sqrt{\alpha n}$ for some $\alpha
  > 0$, then
  \[ \log \Gamma _{\Semigroup{}} \n[X] \sim \frac {\alpha}{1 - X} \]
  as $X \to 1 ^{-}$.

  \par If $\Delta _{\Semigroup{}} \n \sim \frac{c}{1 - X}$ as $X \to 1 ^{-}$,
  then $\GrowthSemigroup{} \n \sim cn$.
\end{theorem}

\subsection{Growth of Mealy automata and of automatic transformation
  semigroups they define} \label{sect:growth_automata}

\begin{definition} \label{def:transformation_semigroup}
  Let $\Automaton = \MAutomaton{n}{m}{}$ be a Mealy automaton. The
  semigroup
\[\Semigroup{\Automaton} = \mathop{sg} \Argument{ f_{q_0}, f_{q_1}, \dots , f
  _{q_{n - 1}}}\] is called the \textit{semigroup of automatic
  transformations defined by $\Automaton$}.
\end{definition}

For an invertible Mealy automaton, let us examine the group of
transformations it defines.  Let $\Automaton$ be a Mealy automaton,
let $\Semigroup{\Automaton}$ be the semigroup defined by $\Automaton$,
and let us denote the growth function and the spherical growth
function of $\Semigroup{\Automaton}$ by the symbols
$\GrowthSemigroup{\Automaton}$ and
$\SphericalGrowthSemigroup{\Automaton}$, respectively. From
Definition~\ref{def:transformation_semigroup} we have

\begin{proposition}[\cite{MR89f:20065}]\label{prop:equivalent_growth_functions}
  For any $\PositiveInteger$ the value $\GrowthAutomaton{} \n$ equals
  the number of those elements of $\Semigroup{\Automaton}$ that can be
  presented as a product of length $n$ in the generators $\left\{
    f_{q_0}, f_{q_1}, \dots, f_{q_{n - 1}} \right\}$, i.e.
  \[\GrowthAutomaton{} \n = \SphericalGrowthSemigroup{\Automaton}
  \n, \, \PositiveInteger.\]
\end{proposition}

\par
From this proposition and~\eqref{eq:growths_semigroup} it follows that
$\GrowthAutomaton{} \n \le \GrowthSemigroup{\Automaton} \n$ for any
$\PositiveInteger$.

\begin{proposition}[\cite{reznykov:phd}] \label{prop:similar_automata}
  Let $\Automaton[i] = \MAutomaton{n}{m}{i}$ for $i = 1, 2$ be two
  similar Mealy automata. Then these automata define isomorphic
  automatic transformation semigroups and have the same growth
  function.
\end{proposition}

\subsection{Hausdorff dimension} \label{subsect:Hausdorff_dimension}

\par We introduce now the Hausdorff dimension of semigroups acting on
trees. This topic was already extensively studied for
groups~\cites{abercrombie:subgroups,barnea-s:hausdorff}.

\par Let $\Semigroup{}$ be a semigroup acting on a tree $\FiniteSet$. This
action extends to an action on the boundary $\InfiniteSet$ of the tree. This
space has the topology of the Cantor set, and can be given the natural metric
\[
    d\Pair{v}{w} = \sup \left\{ m ^{-n} : \, v_n \neq w_n \right\},
\] where $v = v_0 v_1 v_2 \dots, w = w_0 w_1 w_2 \dots \in \InfiniteSet$.
This metric induces the Cantor topology on $\InfiniteSet$, and turns it into a
compact space of diameter $1$.

\par The semigroup $\Semigroup{}$ is a subset of the semigroup of tree
endomorphisms of $\FiniteSet$, and $\EndFinite$ has the natural function
(compact-open) topology. The natural metric on $\EndFinite$ is
\[
    d\Pair{g}{h} = \sup \left\{ \Card{\EndLevels{m}{n}} ^{-1} : \, \text{ there
    exists } v \in \Alphabet{m} ^n \text{ with } v ^g \neq v ^h \right\},
\] where $g, h$ are arbitrary homomorphisms and $\Alphabet{m} ^n$ denotes the
first $n$ levels of the tree $\FiniteSet$. This induces on $\EndFinite$, and
therefore on $\Semigroup{}$, the Cantor topology, and turns $\EndFinite$ into a
compact space of diameter $1$.

\par Furthermore, $\EndFinite$ has Hausdorff dimension $1$, since it
is covered by $\Card{\EndLevels{m}{n}}$ subspaces of diameter
$\Card{\EndLevels{m}{n}}^{-1}$.  Let $\Quotient$ denote the image of
$\Semigroup{}$ in $\EndLevels{m}{n}$; then we define the Hausdorff dimension of
$\Semigroup{}$ as
\[
    H\!\dim{\Semigroup{}} = \liminf_{n \to \infty} \frac{\log
    \Card{\Quotient}}{\log \Card{\EndLevels{m}{n}}}.
\] This is a number in the interval $\left[ 0, 1 \right]$ which measures the
proportion of $\EndFinite$ occupied by $\Semigroup{}$.

\par Let us compute $\Card{\EndLevels{m}{n}}$: such an endomorphism is
determined by an endomorphism of $\Alphabet{m}$ (there are $m^m$ of them), and
$m$ endomorphisms in $\EndLevels{m}{n - 1}$; we arrive at the recursive formula
\[
    \Card{\EndLevels{m}{n}} = m^m \Card{\EndLevels{m}{n - 1}}^m = m^{m
    \frac{m^n - 1}{m - 1}}.
\]

\section{The semigroup $\Semigroup{\Intermediate{2}}$}

\subsection{Properties of automatic transformations}

\par For $i=0,1$ let us denote the automatic transformation
$f_{q_i, \Intermediate{2}}$ by the symbol $f_i$. The unrolled forms of
the automatic transformations $f_0$ and $f_1$ are the following:
\begin{equation} \label{eq:unrolled_form}
\begin{aligned}
f_0 & = \TF{0}{0} \TX{1}{0}, & f_1 & = \TF{1}{0} \TX{1}{1}.
\end{aligned}
\end{equation}
From~\eqref{eq:unrolled_form} the following equalities hold:
\begin{equation} \label{eq:unrolled_form_2}
\begin{aligned}
f_0 ^2 &= \Pair{f_0 ^2}{f_0 ^2} \TX{0}{1}, & f_0 f_1 &= \Pair{f_0 f_1}{f_0 ^2}
\TX{0}{0},\\
f_1 ^2 &= \Pair{f_0 f_1}{f_0 ^2} \TX{1}{1}, & f_1 f_0 &= \Pair{f_0 ^2}{f_1 f_0}
\TX{1}{1};
\end{aligned}
\end{equation}
whence we have
\begin{lemma} \label{lem:involution}
  The automatic transformation $f_0$ is an involution.
\end{lemma}

\par From Lemma~\ref{lem:involution} and~\eqref{eq:unrolled_form_2},
the following equalities hold for any $p \ge 1$:
\begin{subequations} \label{eq:unrolled_form_p}
\begin{align}
  \label{eq:unrolled_form_p_01} \left( f_0 f_1 \right) ^p &= \Pair{\left( f_0
  f_1 \right) ^p}{\left( f_0 f_1 \right) ^{p - 1}} \TX{0}{0}, \\
  \label{eq:unrolled_form_p_10} \left( f_1 f_0 \right) ^p &= \Pair{\left( f_1
  f_0 \right) ^{p - 1}}{\left( f_1 f_0 \right) ^p} \TX{1}{1}, \\
  \label{eq:unrolled_form_p_011} \left( f_0 f_1 \right) ^p f_1 &= \Pair{\left(
  f_0 f_1 \right) ^{p - 1} f_1}{\left( f_0 f_1 \right) ^{p - 1} f_0} \TX{0}{0}.
\end{align}
Here we assume $f ^0 = 1$ for an arbitrary automatic transformation $f$.
\end{subequations}

\begin{lemma} \label{lem:defining_relations}
  In the semigroup $\Semigroup{\Intermediate{2}}$ the following
  relations hold:
\begin{equation} \label{eq:defining_relations}
r_p: \, f_1 \left( {f_0 f_1} \right) ^p \left( {f_1 f_0} \right) ^p f_1 ^2  =
f_1 \left( {f_0 f_1 } \right) ^p \left( {f_1 f_0 } \right) ^p,
\end{equation}
for all $p \ge 0$.
\end{lemma}

\begin{proof} Let us prove the lemma by induction on $p$. For $p = 0$
from~\eqref{eq:unrolled_form_2} follows \[f_1 ^3 = \Pair{f_0 f_1}{1} \TX{1}{1}
\cdot \Pair{f_1}{f_0} \TX{1}{1} = \Pair{f_1}{f_0} \TX{1}{1} = f_1.\] For $p >
1$ from~\eqref{eq:unrolled_form_2},~\eqref{eq:unrolled_form_p_01}
and~\eqref{eq:unrolled_form_p_10} we have
\begin{multline*}
f_1 \left( f_0 f_1 \right) ^p \left( f_1 f_0 \right) ^p = \\ = \Pair{f_1}{f_0}
\TX{1}{1} \cdot \Pair{\left( f_0 f_1 \right) ^p}{\left( f_0 f_1 \right) ^{p -
1}} \TX{0}{0} \cdot \Pair{\left( f_1 f_0 \right) ^{p - 1}}{\left( f_1 f_0
\right) ^p} \TX{1}{1} = \\ = \Pair{f_1 \left( f_0 f_1 \right) ^{p - 1} \left(
f_1 f_0 \right) ^{p - 1}}{f_1 \left( f_0 f_1 \right) ^{p - 1} \left( f_1 f_0
\right) ^p} \TX{1}{1},
\end{multline*}
and
\begin{multline*}
f_1 \left( f_0 f_1 \right) ^p \left( f_1 f_0 \right) ^p f_1 ^2 = \Pair{f_1
\left( f_0 f_1 \right) ^{p - 1} \left( f_1 f_0 \right) ^{p} f_0 f_1}{f_1 \left(
f_0 f_1 \right) ^{p - 1} \left( f_1 f_0 \right) ^p} \TX{1}{1} = \\ = \Pair{f_1
\left( f_0 f_1 \right) ^{p - 1} \left( f_1 f_0 \right) ^{p - 1} f_1 ^2}{f_1
\left( f_0 f_1 \right) ^{p - 1} \left( f_1 f_0 \right) ^p} \TX{1}{1}.
\end{multline*}
By the induction hypothesis, the right-hand sides of both equalities
define the same automatic transformation, so the lemma holds.
\end{proof}

\begin{remark} \label{rem:even_reducing}
  Application of any defining relation to an arbitrary semigroup word
  changes the length of this word by an even number.
\end{remark}

\begin{remark}
  The relation $r_p$ for all $p \ge 1$ can be written in the following way
\begin{equation*} \label{eq:defining_relations2}
    r_p: \, f_1 \left( {f_0 f_1} \right) ^p f_1 \left( {f_0 f_1} \right) ^p f_1
    = f_1 \left( {f_0 f_1 } \right) ^p f_1 \left( {f_0 f_1 } \right) ^{p - 1}
    f_0.
\end{equation*}
In sequel, we will use both presentations of the relations $r_p$.
\end{remark}

\begin{lemma} \label{lem:left_side_zero}
  For any $\PositiveInteger$ the element $f_1 \left( f_0 f_1 \right)
  ^{n - 1}$ is a left-side zero in the semigroup $\Quotient$. That is,
  the relations
  \begin{align} \label{eq:zero_relations}
    f_1 \left( f_0 f_1 \right) ^{n - 1} f_0 &= f_1 \left( f_0 f_1 \right) ^{n -
    1}, & f_1 \left( f_0 f_1 \right) ^{n - 1} f_1 &= f_1 \left( f_0 f_1 \right)
    ^{n - 1},
  \end{align}
  hold in the semigroup $\Quotient$.
\end{lemma}

\begin{proof}
  It is enough to show that the image of an arbitrary word $u \in \Alphabet{2}
 ^n$ under the action of $f_1 \left( f_0 f_1 \right) ^{n - 1}$ does not depend
 on $u$. Indeed, from~\eqref{eq:unrolled_form_p_01} for any $p > 0$ follows
\[  f_1 \left( f_0 f_1 \right) ^p = \Pair{f_1 \left( f_0 f_1 \right) ^p}{f_1
    \left( f_0 f_1 \right) ^{p - 1}} \TX{1}{1}.
\] Let us write the word $u$ as \[
    u = x_0 ^{t_1} x_1 ^{t_2} x_0 ^{t_3} x_1 ^{t_4} \dots x_0 ^{t_{2k - 1}}
    x_1 ^{t_{2k}},
 \] where $k > 0$, $t_1, t_{2k} \ge 0$, $t_i > 0$, $2 \le i \le 2k -
1$, $\sum \limits _{i = 1} ^{2k} t_i = n$. For $u = x_1 ^n$ we have: \[
    f_1 \left( f_0 f_1 \right) ^{n - 1} \n[u] = x_1 ^{n - 1} \cdot f_1 \n[x_1]
    = x_1 ^n.
\] Otherwise, $\sum \limits _{i = 1} ^{k} {t_{2i}} < n$ and the equalities
hold:
\begin{multline*}
    f_1 \left( f_0 f_1 \right) ^{n - 1} \Argument{u} = x_1 ^{t_1} \cdot f_1 \left(
    f_0 f_1 \right) ^{n - 1} \Argument{x_1 ^{t_2} x_0 ^{t_3} x_1 ^{t_4} \dots x_0
    ^{t_{2k - 1}} x_1 ^{t_{2k}}} = \\
    = x_1 ^{t_1 + t_2} \cdot f_1 \left( f_0 f_1 \right) ^{n - 1 - t_2} \Argument{
    x_0 ^{t_3} x_1 ^{t_4} \dots x_0 ^{t_{2k - 1}} x_1 ^{t_{2k}}} = \\
    = \dots = x_1 ^{t_1 + t_2 + t_3 + \dots + t_{2k - 1}} \cdot f_1 \left(
    f_0 f_1 \right) ^{n - 1 - t_2 - t_4 - \dots - t_{2k - 2}} \Argument{x_1 ^{t_{2k}}}
    = x_1 ^n.
\end{multline*}
Therefore $f_1 \left( f_0 f_1 \right) ^{n - 1} \n[u] = x_1 ^n$ and
the lemma holds.
\end{proof}

\subsection{Normal forms}

\begin{proposition}[\cite{reznykov-s:interm2x2}] \label{prop:normal_form}
  Every $\s \in \Semigroup{\Intermediate{2}}$ admits a unique
  minimal-length representation as a word of the form $1$, $f_0$, or
\begin{equation} \label{eq:normal_form}
f_0 ^{\epsilon _1} f_1 \Argument{f_0 f_1} ^{p_1} f_1 \Argument{f_0 f_1} ^{p_2}
f_1 \dots \Argument{f_0 f_1} ^{p_k} f_1 \Argument{f_0 f_1} ^{p_{k + 1}} f_0
^{\epsilon _2},
\end{equation}
where $\epsilon _1, \epsilon _2 \in \left\{ 0,1 \right\}$, $k \ge 0$,
$0 \le p_1 < p_2 < \dots < p_k$, and $p _{k + 1} \ge 0$.
\end{proposition}

\begin{proof}
Let $\s \in \Semigroup{\Intermediate{2}}$ be an arbitrary semigroup element,
written in the following way: \[
    f_0 ^{p_0} f_1 ^{p_1} f_0 ^{p_2} f_1 ^{p_3} \dots f_0 ^{p_{2k}} f_1
    ^{p_{2k + 1}},
\] where $k \ge 0$, $p_0 \ge 0$, $p_{2k + 1} \ge 0$, $p_i > 0$, $i = 1, 2,
\dots, 2k$. The relation $f_0 ^2 = 1$ implies that there can never be two
consecutive $f_0$'s in a reduced word, and the relation $r_0$ is $f_1 ^3 =
f_1$, so there can never be three consecutive $f_1$'s.

\par If the representation of  $\s$ contains at least one symbol
$f_1$, then it can be written in the form
\begin{equation} \label{eq:almost_normal_form}
    \s = f_0 ^{\epsilon _1} f_1 \Argument{f_0 f_1} ^{p_1} f_1 \Argument{f_0 f_1}
    ^{p_2} f_1 \dots \Argument{f_0 f_1} ^{p_k} f_1 \Argument{f_0 f_1} ^{p_{k +
    1}} f_0 ^{\epsilon _2},
\end{equation}
where $\epsilon _1, \epsilon _2 \in \left\{ 0,1 \right\}$, $k \ge 0$, $p_1,
p_{k + 1} \ge 0$, $p_i > 0$, $2 \le i \le k$. Furthermore if $p_i \ge p _{i +
1}$ for some $i \in \left\{ 1, 2, \dots, k - 1 \right\}$, we have the relation
\[
    r _{p_{i + 1}} : \, f_1 \Argument{f_0 f_1} ^{p_{i + 1}} f_1 \Argument{f_0
    f_1} ^{p_{i + 1}} f_1 = f_1 \Argument{f_0 f_1} ^{p_{i + 1}} f_1
    \Argument{f_0 f_1} ^{p_{i + 1} - 1} f_0,
\] and therefore the representation can be shortened. Then the semigroup word
$\s$ is irreducible if and only if for all $i = 1, 2, \dots, {k - 1}$
the inequality $p_i < p_{i + 1}$ holds, that is $0 \le p_1 < p_2 <
\dots < p_k$.
\end{proof}

\par In~\cite{reznykov-s:interm2x2} the algorithm of reducing an arbitrary
semigroup word to normal form is considered. Let $\s$ be an arbitrary
semigroup word over the alphabet $\left\{ f_0, f_1 \right\}$. It can
be reduced to normal form by the following steps:

\begin{enumerate}
\item \label{step:first} The word $\s$ is reduced by the defining
  relation $f_0 ^2 = 1$;

\item \label{step:second} The word $\s$ is reduced by the defining
  relation $r_0$;

\item After steps~\ref{step:first} and~\ref{step:second} the word is
  written as~\eqref{eq:almost_normal_form};
  
\item If for all $i = 1, 2, \dots, {k - 1}$ the numbers $p_i$
  in~\eqref{eq:almost_normal_form} satisfy the inequalities $p_i <
  p_{i + 1}$, then the algorithm finishes, otherwise it goes to the
  next step;
  
\item For the first pair of exponents $p _j$ and $p _{j + 1}$, with $1
  \leq j \leq {k - 1}$, such that $p _j \geq p _{j + 1}$, the subword
  $f_1 ^2$ of length $2$ is canceled in $\s$, by the application of
  the relation $r _{p _{j + 1}}$;

\item \label{step:last} Go to step~\ref{step:first}.
\end{enumerate}

\begin{proposition}[\cite{reznykov-s:interm2x2}]
  \label{prop:reducing_algorithm}
  The algorithm with steps~\ref{step:first}--\ref{step:last} reduces
  an arbitrary semigroup word $\s$ to its normal form in no more than
  $\left[ \frac{\left| \s \right|}{2} \right]$ steps.
\end{proposition}

\begin{lemma} \label{prop:normal_form_quotient}
  For any $n \ge 1$ an arbitrary element $\s$ of $\Quotient$ equals
  $1$, $f_0$, or can be written in normal form
  \begin{equation} \label{eq:normal_form_quotient}
    f_0 ^{\epsilon _1} f_1 \Argument{f_0 f_1} ^{p_1} f_1 \Argument{f_0 f_1}
    ^{p_2} f_1 \dots \Argument{f_0 f_1} ^{p_k} f_1 \Argument{f_0 f_1} ^{p_{k +
    1}} f_0 ^{\epsilon _2},
  \end{equation}
where $\epsilon _1, \epsilon _2 \in \left\{ 0,1 \right\}$, $0 \le k$, $0 \le
p_1 < p_2 < \dots < p_k < {n - 1}$, and $0 \le p _{k + 1} + \epsilon_2 \le {n
- 1}$.
\end{lemma}

\begin{proof}
Let us fix a number $n \ge 1$. Let $\s$ be an arbitrary word of normal
form~\eqref{eq:normal_form}: \[
    \s = f_0 ^{\epsilon _1} f_1 \Argument{f_0 f_1} ^{p_1} f_1 \Argument{f_0
    f_1} ^{p_2} f_1 \dots \Argument{f_0 f_1} ^{p_k} f_1 \Argument{f_0 f_1}
    ^{p_{k + 1}} f_0 ^{\epsilon _2},
\] where $\epsilon _1, \epsilon _2 \in \left\{ 0,1 \right\}$, $0 \le k$, $0 \le
p_1 < p_2 < \dots < p_k$, and $0 \le p _{k + 1}$. If $p_{i} \ge {n - 1}$ for
some $i$, then the semigroup word may be shortened by using the
relations~\eqref{eq:zero_relations}:
\begin{multline*}
    \s = f_0 ^{\epsilon _1} f_1 \Argument{f_0 f_1} ^{p_1} f_1 \Argument{f_0
    f_1} ^{p_2} f_1 \dots \Argument{f_0 f_1} ^{p_i} f_1 \dots \Argument{f_0
    f_1} ^{p_k} f_1 \Argument{f_0 f_1} ^{p_{k + 1}} f_0 ^{\epsilon _2} =\\
    = f_0 ^{\epsilon _1} f_1 \Argument{f_0 f_1} ^{p_1} f_1 \dots \Argument{f_0 f_1}
    ^{p_{i - 1}} f_1 \Argument{f_0 f_1} ^{n - 1} \Argument{f_0 f_1} ^{p_i - n
    + 1} f_1 \dots \Argument{f_0 f_1} ^{p_k}
    f_1 \Argument{f_0 f_1} ^{p_{k + 1}} f_0 ^{\epsilon _2}= \\
    = f_0 ^{\epsilon _1} f_1 \Argument{f_0 f_1} ^{p_1} f_1 \Argument{f_0 f_1}
    ^{p_2} f_1 \dots \Argument{f_0 f_1} ^{p_{i - 1}} f_1 \Argument{f_0 f_1}
    ^{n - 1}.
\end{multline*}
This gives the requirements $0 \le p_1 < p_2 < \dots < p_k < {n -
  1}$.  Similarly, the end of $\s$, the subword $\Argument{f_0 f_1}
^{p_{k + 1}} f_0 ^{\epsilon _2}$, should be no longer than
$\Argument{f_0 f_1} ^{n - 1}$. Hence, the requirement $p _{k + 1} +
\epsilon_2 \le {n - 1}$ should be satisfied.
\end{proof}

\subsection{Proof of Theorem~\ref{th:semigroup}}


\par Let $k \in \Natural$ be an arbitrary positive integer, and let us denote
its remainder modulo $2$ by the symbol $\Remainder{k}$.

\begin{proposition} \label{prop:test_word}
  Let $\s \in \Semigroup{\Intermediate{2}}$ be an arbitrary element
  such that \[ \s = f_0 f_1 \Argument{f_0 f_1} ^{p_1} f_1
  \Argument{f_0 f_1} ^{p_2} f_1 \dots \Argument{f_0 f_1} ^{p_k} f_1 ,
  \] where $k \ge 1$,
  $0 \le p_1 < p_2 < \dots < p_k$. Then \[
    \s \Argument{x_0 ^\ast} = x_0 ^{p_1 + 1} x_1 ^{p_2 - p_1} x_0 ^{p_3 - p_2}
    \dots x_{1 - \Remainder{k}} ^{p_k - p_{k - 1}} x_{\Remainder{k}} ^\ast.
  \]
\end{proposition}

\begin{proof}
  Let $u \in \InfiniteSet[2]$ be an arbitrary word, and let $t_2 \ge
  t_1 \ge 0$ be arbitrary integers. Then
  from~\eqref{eq:unrolled_form_p_011} we have
\begin{multline*}
  \left( f_0 f_1 \right) ^{t_1} f_1 \Argument{x_0 ^{t_2} x_1 u} = x_0
  \cdot \left( f_0 f_1 \right) ^{t_1 - 1} f_1 \Argument{x_0 ^{t_2 - 1}
    x_1 u} = \dots = \\ = x_0 ^{t_1} \cdot f_1 \Argument{x_0 ^{t_2 -
      t_1} x_1 u} = x_0 ^{t_1} x_1 ^{t_2 - t_1 + 1} \cdot f_0
  \Argument{u}.
\end{multline*}

\par Let us prove the lemma by induction on $k$. For $k = 1$
from~\eqref{eq:unrolled_form_p_011} follows \[
    \left( f_0 f_1 \right) ^{p_1 + 1} f_1 \Argument{x_0 ^\ast} = x_0 ^{p_1 + 1}
    \cdot f_1 \Argument{x_0 ^\ast} = x_0 ^{p_1 + 1} x_{1} ^\ast.
\] For $k > 1$ we have
\begin{align*}
  \s \Argument{x_0 ^\ast} & = f_0 f_1 \Argument{f_0 f_1} ^{p_1} f_1 \Argument{
  f_0 f_1 \Argument{f_0 f_1} ^{p_2 - 1} f_1 \Argument{f_0 f_1} ^{p_3} f_1
  \dots \Argument{f_0 f_1} ^{p_k} f_1 \Argument{x_0 ^\ast}} = \\
  & = f_0 f_1 \Argument{f_0 f_1} ^{p_1} f_1 \Argument{x_0
  ^{p_2} x_1 ^{p_3 - p_2 + 1} x_0 ^{p_4 - p_3} \dots x_{1 - \Remainder{k - 1}}
  ^{p_k - p_{k - 1}} x_{\Remainder{k - 1}} ^\ast} = \\
  & = x_0 ^{p_1 + 1} x_1 ^{p_2 - p_1} \cdot f_0 \Argument{ x_1 ^{p_3 - p_2} x_0
  ^{p_4 - p_3} \dots x_{\Remainder{k}} ^{p_k - p_{k - 1}} x_{1 -
  \Remainder{k}} ^\ast} = \\
  & = x_0 ^{p_1 + 1} x_1 ^{p_2 - p_1} x_0 ^{p_3 - p_2}
  \dots x_{1 - \Remainder{k}} ^{p_k - p_{k - 1}} x_{\Remainder{k}} ^\ast,
\end{align*}
and the lemma holds.
\end{proof}

\begin{corollary} \label{cor:test_word_quotient}
Let $\PositiveInteger$ be any, and let $\s$ be an semigroup element, written in
the following form~\eqref{eq:normal_form_quotient}: \[
    \s = f_0 f_1 \Argument{f_0 f_1} ^{p_1} f_1 \Argument{f_0 f_1} ^{p_2} f_1
    \dots \Argument{f_0 f_1} ^{p_k} f_1 \Argument{f_0 f_1} ^{p_{k + 1}},
\] where $0 \le k$, $0 \le p_1 < p_2 < \dots < p_k < {n - 1}$, and $0 \le p
_{k + 1} \le {n - 1}$. Then \[
    \s \Argument{x_0 ^n} = x_0 ^{p_1 + 1} x_1 ^{p_2 - p_1} x_0 ^{p_3 - p_2}
    \dots x_{1 - \Remainder{k}} ^{p_k - p_{k - 1}} x_{\Remainder{k}} ^{n - p_k
    - 1}.
\]
\end{corollary}

\begin{proof}
Let us fix an integer $n \ge 1$. From~\eqref{eq:unrolled_form_p_01} for any $p
\ge 0$ we have
 \[
 \Argument{f_0 f_1} ^p \Argument{x_0 ^n} = x_0 ^n.
\] Therefore, for $k = 0$ we have \[
    \s \Argument{x_0 ^n} = f_0 f_1 \Argument{f_0 f_1} ^{p_{k + 1}}
    \Argument{x_0 ^n} = x_0 ^n,
\] and when $k > 0$ from Proposition~\ref{prop:test_word} we have
\begin{align*}
    \s \Argument{x_0 ^n} & = f_0 f_1 \Argument{f_0 f_1} ^{p_1} f_1
    \Argument{f_0 f_1} ^{p_2} f_1 \dots \Argument{f_0 f_1} ^{p_k} f_1
    \Argument{x_0 ^n} = \\
    & = x_0 ^{p_1 + 1} x_1 ^{p_2 - p_1} x_0 ^{p_3 - p_2} \dots x_{1 -
    \Remainder{k}} ^{p_k - p_{k - 1}} x_{\Remainder{k}} ^{n - \Argument{p_1 +
    1} - \sum \limits _{i = 2} ^{k} {\Argument{p_i - p_{i - 1}}}} = \\
    & = x_0 ^{p_1 + 1} x_1 ^{p_2 - p_1} x_0 ^{p_3 - p_2} \dots x_{1 -
    \Remainder{k}} ^{p_k - p_{k - 1}} x_{\Remainder{k}} ^{n - p_k - 1}.
    \qedhere
\end{align*}
\end{proof}

\begin{proposition} \label{prop:infinite_presented}
  The infinite system of relations \[
    f_0 ^2 = 1, r_0 ,r_1 ,r_2 , \dots,
  \] is minimal, that is none of the relations follows from the
  others.
\end{proposition}

\begin{proof}
  Let us show that the relation \[f_0 ^2 = 1\] does not follow from
  the relations $\left\{ r_p, p \ge 0 \right\}$. Indeed, each relation
  $r_p$, for $p \ge 0$, includes the symbol $f_1$ in both its left-
  and right-hand side, and therefore it cannot be applied to $f_0 ^2 =
  1$.

\par Moreover, the relation \[r_0 :\, f_1 ^3 = f_1 \] does not follow from the
set of relations $f_0 ^2 = 1, \left\{ r_p, p \ge 1 \right\}$, either.
Let us consider its right-hand side, the semigroup word $f_1$. The
unique relation which may be applied to it is $f_0 ^2 = 1$; and the
set of semigroup words equivalent to $f_1$ is described in the
following way:
\[ f_0 ^{2p_1} f_1 f_0 ^{2p_2}, \, p_1, p_2 \ge 0.
\] Obviously, this set does not include the semigroup word $f_1 ^3$.

\par Let us denote the left- and right-hand sides of the relation $r_p$,
for $p > 0$, by the symbols $w_p$ and $v_p$ respectively, that is
\begin{align*}
    w_p & = f_1 \left( f_0 f_1 \right) ^p f_1 \left( f_0 f_1 \right) ^{p}
    f_1,\\
    v_p & = f_1 \left( f_0 f_1 \right) ^p f_1 \left( f_0 f_1 \right) ^{p - 1}
    f_0.
\end{align*}
Let us fix a positive integer $\ell \ge 1$ and prove that the set of
semigroup words equivalent to $v_\ell$, obtained by applying the
relations $f_0 ^2 = 1, \left\{ r_p, p \ge 0, p \neq l \right\}$, does
not include any semigroup words which end in the symbol $f_1$.


\par Let us consider the set of semigroup words
\begin{equation*}
    \Omega _i = \left\{ {\begin{array}{*{20}c}
        {f_1 f_0 ^{1 + 2t_1} f_1 f_0 ^{1 + 2t_2} \dots f_1 f_0 ^{1 + 2t_{i - 1}}
        f_1 f_0 ^{2t_{i}}} & \vline & {t_1, t_2, \dots, t_{i} \ge 0}\\
    \end{array}} \right\},
\end{equation*}
where $i > 0$. All words in the set $\Omega _i$, for $i > 0$, are
pairwise equivalent, and let us choose the word of minimal length \[
\omega _i = f_1 f_0 f_1 f_0 \dots f_1 f_0 f_1 = f_1 \Argument{f_0 f_1}
^{i - 1}
\] as the representative of $\Omega_i$. For $i = 0$ let us consider
the set of words
\begin{align*}
    \Omega _0 = \left\{ {\begin{array}{*{20}c}
        {f_0 ^{2t_1}} & \vline & {t_1 \ge 0}\\
    \end{array}} \right\},
\end{align*}
with representative $\omega _0 = 1$.

\par Let $\s \in \Semigroup{\Intermediate{2}}$ be an arbitrary semigroup
element. It can be ambiguously written in the following way \[
    f_0 ^{\epsilon_1} \nu _{0} \nu _{1} \dots \nu _{k} f_0
    ^{\epsilon_2},
\] where $k \ge 0$, $\epsilon_1, \epsilon_2 \in \left\{ 0, 1 \right\}$, $\nu _j
\in \Omega _{i_j}$, $j = 0, 1, \dots, {k}$, $i_0 = 0$, $i_j > 0$, $j = 1, 2,
\dots, {k}$, and if $k = 0$ let $\epsilon_1 \le \epsilon_2$, $\epsilon_1 +
\epsilon_2 \le 1$. Using only the relation $f_0 ^2 = 1$, the element $\s$ can
be ambiguously reduced to the following product
\begin{equation*} \label{eq:word_omega_product}
    \s = f_0 ^{\epsilon_1} \omega _{i_0} \omega _{i_1} \dots \omega _{i_{k}}
    f_0 ^{\epsilon_2},
\end{equation*}
where requirements on parameters are listed above.

\par Let us consider the set \[
    \Upsilon \n[\s] = \left\{ {\begin{array}{*{20}c}
        {\sum \limits _{j = 0} ^{l} {\Argument{-1} ^{j + 1} i_j}} & \vline & {l
        = 0, 1, 2, \dots, k}\\
    \end{array}} \right\}.
\] \emph{The width} of the semigroup word $\s$ is the positive integer \[
    w \n[\s] = \max \Upsilon \n[\s] - \min \Upsilon \n[\s].
\] Let us note that $\s$ has width $0$ if and only if $\s = f_0 ^p$ for some $p \ge
0$.

\par The relations $r_p$, for $p = 0, 1, \dots$, have the following
representations:
\begin{equation} \label{eq:relations_omega_product}
\begin{aligned}
    r_0 & : \, \omega _{0} \omega _{1} \omega _{1} \omega _{1} = \omega _{0}
    \omega _{1};\\
    r_p & : \, \omega _{0} \omega _{p + 1} \omega_{p + 1} \omega_{1} = \omega
    _{0} \omega_{p + 1} \omega_{p} f_0, \quad p > 0.
\end{aligned}
\end{equation}
Obviously, the left- and right-hand sides of $r_p$ have the same width
$\Argument{p + 1}$, for all $p > 0$. Moreover, both sides of the
relation $f_0 ^2 = 1$ have the same width $0$, too.

\par From~\eqref{eq:relations_omega_product} it follows that the application of
relations $f_0 ^2 = 1$ or $r_p$, for $p = 0, 1, \dots$ does not change
the width of $\s$, and the relation $r_p$ can be applied to $\s$ if
and only if $0 \le p \le w \n[\s] - 1$. Hence, only the relations
\begin{align} \label{eq:used_relations}
    f_0 ^2 & = 1, r_0, r_1, \dots, r_{\ell - 1}
\end{align}
can be applied to the word $v_\ell$ and the words equivalent to it.

\par Let us separate $v_\ell$ into two parts
\begin{align*}
    v_\ell ^{\n[1]} & = \omega_0 \omega _{\ell + 1}, & v_\ell ^{\n[2]} & = \omega_0
    \omega _{l} f_0,
\end{align*}
where $v_\ell = v_\ell ^{\n[1]} \cdot v_\ell ^{\n[2]}$. In addition,
\begin{align*}
    w \n[v_\ell] = w \n[v_\ell ^{\n[1]}] & = \ell + 1, & w \n[v_\ell ^{\n[2]}] & = \ell.
\end{align*}
From Proposition~\ref{prop:normal_form} it follows that all words
$v_\ell$, $v_\ell ^{\n[1]}$ and $v_\ell ^{\n[2]}$ have normal
form~\eqref{eq:normal_form}. If the relation $r_p$
from~\eqref{eq:used_relations} is applied to $v_\ell$, then there are
three possible cases:
\begin{itemize}
\item $w_p$ and $v_p$ belong to $v_\ell ^{\n[1]}$; \item $w_p$ and $v_p$ belong to
$v_\ell ^{\n[2]}$; \item $\omega _0 \omega _{p + 1}$ belongs to $v_\ell ^{\n[1]}$,
and $\omega _0 \omega _{p + 1} \omega_1$ and $\omega _0 \omega _{p} f_0$ belong
to $v_\ell ^{\n[2]}$.
\end{itemize}
As $p < \ell$, the application of a relation
from~\eqref{eq:used_relations} does not change the width of the parts
$v_\ell ^{\n[1]}$ and $v_\ell ^{\n[2]}$. Hence, if $\s$ is an arbitrary word
which is obtained from $v_\ell$ by relations~\eqref{eq:used_relations},
it can be separated into two parts $\s ^{\n[1]}$ and $\s ^{\n[2]}$,
$\s = \s ^{\n[1]} \cdot \s ^{\n[2]}$, where $w \n[\s ^{\n[i]}] = w
\n[v_\ell ^{\n[i]}]$ for $i = 1, 2$.  As $w \n[\s ^{\n[2]}] = \ell$ and the
parities of the number of occurrences of $f_0$ in $\s ^{\n[2]}$ and
$v_\ell ^{\n[2]}$ coincide, $\s ^{\n[2]}$ ends on the symbol $f_0$.
Therefore the word $\s = \s ^{\n[1]} \cdot \s ^{\n[2]}$ ends in
$f_0$ too, and the word $\omega _{0} \omega _{p + 1} \omega_{p + 1}
\omega_{1}$, which ends on $f_1$, is not equivalent to $v_\ell$.
\end{proof}

\begin{proof}[Proof of Theorem~\ref{th:semigroup}] From
  Lemmas~\ref{lem:involution} and~\ref{lem:defining_relations} it
  follows that in the semigroup $\Semigroup{\Intermediate{2}}$ the
  relations $f_0 ^2 = 1$ and $r_p$, for $p \ge 0$, hold. In
  Proposition~\ref{prop:normal_form} it is proved that, using these
  relations, each element can be unambiguously reduced to normal form.
  It is enough to show that semigroup elements, which are written in
  different normal forms, define different automatic transformations
  over the set $\InfiniteSet[2]$.

\par Let  $\s _1,\s _2$ be arbitrary semigroup elements, written in normal
form. As $f_0$ is a bijection and $f_1$ is not a bijection, then any
semigroup word which includes the symbol $f_1$ defines an automatic
transformation which is not a bijection, and therefore differs from
both transformations $1$ and $f_0$.  Due to this remark, it is enough
to consider elements in normal form~\eqref{eq:normal_form}. Let us
write
\begin{align*}
    \s_1 & = f_0 ^{\epsilon _1} f_1 \Argument{f_0 f_1} ^{p_1} f_1 \Argument{f_0
    f_1} ^{p_2} f_1 \dots \Argument{f_0 f_1} ^{p_k} f_1 \Argument{f_0 f_1}
    ^{p_{k + 1}} f_0 ^{\epsilon _2},\\
    \s_2 & = f_0 ^{\mu _1} f_1 \Argument{f_0 f_1} ^{t_1} f_1 \Argument{f_0 f_1}
    ^{t_2} f_1 \dots \Argument{f_0 f_1} ^{t_\ell} f_1 \Argument{f_0 f_1} ^{t_{\ell +
    1}} f_0 ^{\mu _2},
\end{align*}
where $\epsilon _1, \epsilon _2, \mu _1, \mu _2 \in \left\{ 0,1
\right\}$, $k, \ell \ge 0$, $0 \le p_1 < p_2 < \dots < p_k$, $0 \le t_1
< t_2 < \dots < t_\ell$, and $p _{k + 1} \ge 0$, $t _{\ell + 1} \ge 0$.

\par Let us assume that the elements $\s_1$ and $\s_2$ define the same automatic
transformation over $\InfiniteSet[2]$. Then for any $u \in
\InfiniteSet[2]$ the equality holds
\begin{equation} \label{eq:assumption}
\s_1 \Argument{u} = \s_2 \Argument{u}.
\end{equation}
As $f_0$ is a bijection, the equalities
\begin{align*}
\s_1 &= \s_2, & f_0 \s_1 &= f_0 \s_2, & \s_1 f_0 &= \s_2 f_0,
\end{align*}
hold simultaneously. Moreover, from~\eqref{eq:assumption} for any element $\s_3
\in \Semigroup{\Intermediate{2}}$ it follows that \[
    \s_1 \s_3 \n[u] = \s_2 \s_3 \n[u].
\]

\par Let us consider possible values of $\epsilon _1$ and $\mu _1$.

\begin{enumerate}
\item \label{case:e0m0} $\epsilon _1 = 0$ and $\mu _1 = 0$. Due to the note
  above, this case is equivalent to the case $\epsilon _1 = 1$ and $\mu _1 =
  1$, which is described below.

\item \label{case:e1m0} $\epsilon _1 = 0$ and $\mu _1 = 1$. As $k, l
  \geq 0$, the semigroup words $\s_1$ and $\s_2$ start by the symbols
  $f_1$ and $f_0 f_1$, respectively. For the input word $u = x_1$ we
  have
\begin{align*}
    \s_1 \Argument{x_1} & = f_1 \Argument{ \Argument{f_0 f_1} ^{p_1} f_1 \dots
    \Argument{f_0 f_1} ^{p_{k}} f_1 \Argument{f_0 f_1} ^{p_{k + 1}} f_0
    ^{\epsilon _2} \Argument{x_1}} = x_1,
  \intertext{and}
    \s_2 \Argument{x_1} & = f_0 f_1 \Argument{ \Argument{f_0 f_1} ^{t_1} f_1
    \dots \Argument{f_0 f_1} ^{t_{\ell}} f_1 \Argument{f_0 f_1} ^{t_{\ell + 1}} f_0
    ^{\mu _2} \Argument{x_1}} = f_0 \Argument{x_1} = x_0.
\end{align*}
Therefore, the elements $\s_1$ and $\s_2$ define different automatic
transformations over the set $\InfiniteSet[2]$. The case $\epsilon _1
= 1$ and $\mu _1 = 0$ is similar.

\item \label{case:e1m1} $\varepsilon _1 = 1$ and $\mu _1 = 1$.
\end{enumerate}


Let us assume that $\varepsilon _2 = 0$ and $\mu _2 = 0$.
From~\eqref{eq:assumption} it follows that the elements \[\s_1
\Argument{f_0 f_1} ^\Argument{p_k + t_\ell + 1} f_1 \quad \text{and}
\quad \s_2 \Argument{f_0 f_1} ^\Argument{p_k + t_\ell + 1} f_1\] define
the same automatic transformation. Using
Proposition~\ref{prop:test_word}, we have
\begin{align*}
    \s_1 \Argument{f_0 f_1} ^\Argument{p_k + t_\ell + 1} f_1 \Argument{x_0 ^\ast}
    &= f_0 f_1 \Argument{f_0 f_1} ^{p_1} f_1 \dots \Argument{f_0 f_1} ^{p_k}
    f_1 \Argument{f_0 f_1} ^{p_{k + 1} + p_k + t_\ell + 1} f_1 \Argument{x_0
    ^\ast} = \\
    & = x_0 ^{p_1 + 1} x_1 ^{p_2 - p_1} x_0 ^{p_3 - p_2} \dots x_{1 -
    \Remainder{k}} ^{p_k - p_{k - 1}} x_{\Remainder{k}} ^{p_{k + 1 } + t_\ell + 1}
    x_{1 - \Remainder{k}} ^\ast,\\
    \s_2 \Argument{f_0 f_1} ^\Argument{p_k + t_\ell + 1} f_1 \Argument{x_0 ^\ast}
    &= f_0 f_1 \Argument{f_0 f_1} ^{t_1} f_1 \dots \Argument{f_0 f_1} ^{t_\ell}
    f_1 \Argument{f_0 f_1} ^{t_{\ell + 1} + p_k + t_\ell + 1} f_1 \Argument{x_0
    ^\ast} = \\
    & = x_0 ^{t_1 + 1} x_1 ^{t_2 - t_1} x_0 ^{t_3 - t_2} \dots x_{1 -
    \Remainder{l}} ^{t_\ell - t_{\ell - 1}} x_{\Remainder{l}} ^{t_{\ell + 1 } + p_k + 1}
    x_{1 - \Remainder{l}} ^\ast.
\end{align*}
As the words in the right-hand sides coincide, we obtain the
requirements $k = \ell$, $p_i = t_i$, $i = 1, 2, \dots, k + 1$. This
means that the elements $\s_1$ and $\s_2$ are written in the same
normal form~\eqref{eq:normal_form}.

\par The case $\varepsilon _2 = 1$ and $\mu _2 = 1$ is considered similarly,
because we may consider elements $\s_1 f_0$ and $\s_2 f_0$.

\par Next, let us assume $\epsilon _2 = 0$ and $\mu _2 = 1$ (the case
$\epsilon_2 = 1$ and $\mu _2 = 0$ is considered similarly). If $k =
0$ or $p _{k + 1} > p_k$, then elements
\begin{align*}
    \s_1 f_1 & = f_0 f_1 \Argument{f_0 f_1} ^{p_1} f_1 \Argument{f_0 f_1}
    ^{p_2} f_1 \dots \Argument{f_0 f_1} ^{p_k} f_1 \Argument{f_0 f_1} ^{p_{k +
    1}} f_1,\\
    \s_2 f_1 & = f_0 f_1 \Argument{f_0 f_1} ^{t_1} f_1 \Argument{f_0 f_1}
    ^{t_2} f_1 \dots \Argument{f_0 f_1} ^{t_\ell} f_1 \Argument{f_0 f_1} ^{t_{\ell +
    1}} f_0 f_1,
\end{align*}
are written in normal form~\eqref{eq:normal_form} and define the same
transformation over $\InfiniteSet[2]$. From the proof above in
Case~\ref{case:e1m1}.\ it follows that $k + 1 = \ell$, $p _i = t _i$, $i =
1, 2, \dots, k + 1$, and $t _{\ell + 1} + 1 = 0$; but this contradicts
the condition $t _{\ell + 1} \geq 0$.

\par Let $k > 0$ and $0 \leq p _{k + 1} \leq p _k$. Let us assume $\s_3 = f_1
\Argument{f_0 f_1} ^{p _{k + 1}}$, then the element $\s_2 \s_3$ is
already written in normal form~\eqref{eq:normal_form}:
\begin{align*}
    \s_2 f_1 \Argument{f_0 f_1} ^{p _{k + 1}} & = f_0 f_1 \Argument{f_0 f_1}
    ^{t_1} f_1 \Argument{f_0 f_1} ^{t_2} f_1 \dots \Argument{f_0 f_1} ^{t_\ell}
    f_1 \Argument{f_0 f_1} ^{t_{\ell + 1}} f_0 f_1 \Argument{f_0 f_1} ^{p _{k +
    1}} = \\
    & = f_0 f_1 \Argument{f_0 f_1} ^{t_1} f_1 \Argument{f_0 f_1} ^{t_2} f_1
    \dots \Argument{f_0 f_1} ^{t_\ell} f_1 \Argument{f_0 f_1} ^{t_{\ell + 1} + p_{k
    + 1} + 1}.
\end{align*}


\noindent The element $\s_1 \s_3$ is reduced, and its normal form is the following:
\begin{align*}
    \s_1 f_1 \Argument{f_0 f_1} ^{p _{k + 1}} & = f_0 f_1 \Argument{f_0 f_1}
    ^{p_1} f_1 \Argument{f_0 f_1} ^{p_2} f_1 \dots \Argument{f_0 f_1} ^{p_k}
    f_1 \Argument{f_0 f_1} ^{p_{k + 1}} f_1 \Argument{f_0 f_1} ^{p _{k + 1}} =
    \\
    & = f_0 f_1 \Argument{f_0 f_1} ^{p_1} f_1 \Argument{f_0 f_1} ^{p_2} f_1
    \dots \Argument{f_0 f_1} ^{p_{k - 1}} f_1 \Argument{f_0 f_1} ^{p_k}.
\end{align*}
For the elements $\s_1 \s_3$ and $\s_2 \s_3$, as proved above, the
following requirements hold:
\begin{equation} \label{eq:requirements}
k \ge 1, k - 1 = \ell, p_1 = t_1, p_2 = t_2, \dots, p_{k - 1} = t_\ell, p_k = t_{\ell +
1} + p_{k + 1} + 1.
\end{equation}

\par As $f_0$ is a bijection, a similar reasoning can be carried out for the
elements $\s_1 f_0$ and $\s_2 f_0$, where $\s_1$ and $\s_2$ are
rearranged. For the case $t_{\ell + 1} > t_\ell$ or $\ell = 0$ we obtain a
contradiction with the requirement $p_{k + 1} \geq 0$, and in the case $0
\leq t_{\ell + 1} \leq t_\ell$ and $\ell > 0$ the requirement $\ell - 1 = k$
should be fulfilled, but it contradicts the
requirements~\eqref{eq:requirements}.

\par Thus, the relations $f_0 ^2 = 1, r_0, r_1, \dots$ form the system of
defining relations. In Proposition~\ref{prop:infinite_presented} it is proved
that this system is minimal, and therefore the semigroup
$\Semigroup{\Intermediate{2}}$ is infinitely presented.

\par To solve the word problem in $\Semigroup{\Intermediate{2}}$, it is
necessary to reduce semigroup words $\s_1$ and $\s_2$ to normal
form~\eqref{eq:normal_form}, and then to check them for graphical equality.
From Proposition~\ref{prop:reducing_algorithm} this can be done in no more than
\[\left[ \frac{\left| \s_1 \right|}{2} \right] + \left[ \frac{\left| \s_2
\right|}{2} \right]\] steps, and the word problem is solved in polynomial time.

\par Let us prove the second part of Theorem~\ref{th:semigroup} in a similar
way as the first part. Let us fix the integer $n \ge 1$ and let $\s_1$ and
$\s_2$ are arbitrary elements of the semigroup $\Quotient$. The elements $1$,
$f_0$, $f_1 \dots $ and $f_0 f_1 \dots$ define pairwise distinct
transformations $\TX{0}{1}$, $\TX{1}{0}$, $\TX{1}{1}$, and $\TX{0}{0}$ over the
set $\Alphabet{2} ^1$, respectively. Therefore, using the proof above, it is
enough to consider $n > 1$ and $\s_1,\s_2$ such that they are written in
normal form~\eqref{eq:normal_form_quotient}:
\begin{equation} \label{eq:elements_quotient}
\begin{aligned}
    \s_1 & = f_0 f_1 \Argument{f_0 f_1} ^{p_1} f_1 \Argument{f_0 f_1} ^{p_2}
    f_1 \dots \Argument{f_0 f_1} ^{p_k} f_1 \Argument{f_0 f_1} ^{p_{k + 1}}
    f_0 ^{\epsilon _2},\\
    \s_2 & = f_0 f_1 \Argument{f_0 f_1} ^{t_1} f_1 \Argument{f_0 f_1} ^{t_2}
    f_1 \dots \Argument{f_0 f_1} ^{t_\ell} f_1 \Argument{f_0 f_1} ^{t_{\ell + 1}}
    f_0 ^{\mu _2},
\end{aligned}
\end{equation}
where $\epsilon _2, \mu _2 \in \left\{ 0,1 \right\}$, $k, \ell \ge 0$, $0 \le p_1
< p_2 < \dots < p_k < n - 1$, $0 \le t_1 < t_2 < \dots < t_\ell < n - 1$, and $0
\le p _{k + 1} + \epsilon_2 \le n - 1$, $0 \le t _{\ell + 1} + \mu_2 \le n - 1$.
Let us consider these elements in the same way as it was done for elements of
the semigroup $\Semigroup{\Intermediate{2}}$. Besides, it is enough to consider
the cases $\epsilon _2 = \mu _2 = 0$ and $\epsilon _2 = 1$, $\mu _2 = 0$.

\par Let us assume that $\epsilon _2 = \mu _2 = 0$. Then from
Corollary~\ref{cor:test_word_quotient} for the input word $u = x_0 ^n$
it follows that
\begin{align*}
    \s_1 \Argument{x_0 ^n} & = x_0 ^{p_1 + 1} x_1 ^{p_2 - p_1} x_0 ^{p_3 - p_2}
    \dots x_{1 - \Remainder{k}} ^{p_k - p_{k - 1}} x_{\Remainder{k}} ^{n - p_k
    - 1},\\
    \s_2 \Argument{x_0 ^n} & = x_0 ^{t_1 + 1} x_1 ^{t_2 - t_1} x_0 ^{t_3 - t_2}
    \dots x_{1 - \Remainder{l}} ^{t_\ell - t_{\ell - 1}} x_{\Remainder{l}} ^{n - t_\ell
    - 1},
\end{align*}
and from Assumption~\eqref{eq:assumption} we have the requirements \[
    k = \ell, p_i = t_i, i = 1, 2, \dots, k.
\] With no loss of generality let us assume $0 \le p_{k + 1} < t_{k + 1}$.

\par If $k = 0$ or $p_k < n - 1 - t_{k + 1} + p_{k + 1}$, let us consider the
element $\s_3 = \Argument{f_0 f_1} ^{n - 1 - t_{k + 1}} f_1$. Then
$\s_1 \s_3$ does not reduce, because $p_{k + 1} + n - 1 - t_{k + 1} <
n - 1$, and $\s_2 \s_3$ is reduced to the following element:
\begin{multline*}
    \s_2 \s_3 = f_0 f_1 \Argument{f_0 f_1} ^{p_1} f_1 \Argument{f_0 f_1} ^{p_2}
    f_1 \dots \Argument{f_0 f_1} ^{p_k} f_1 \Argument{f_0 f_1} ^{t_{k + 1}}
    \cdot \Argument{f_0 f_1} ^{n - 1 - t_{k + 1}} f_1 = \\
    = f_0 f_1 \Argument{f_0 f_1} ^{p_1} f_1 \Argument{f_0 f_1} ^{p_2} f_1
    \dots \Argument{f_0 f_1} ^{p_k} f_1 \Argument{f_0 f_1} ^{n - 1}.
\end{multline*}
For the input word $u = x_0 ^n$ we have
\begin{align*}
    \s_1 \s_3 \Argument{x_0 ^n} & = f_0 f_1 \Argument{f_0 f_1} ^{p_1} f_1
    \Argument{f_0 f_1} ^{p_2} f_1 \dots \Argument{f_0 f_1} ^{p_k} f_1
    \Argument{f_0 f_1} ^{n - 1 - t_{k + 1} + p_{k + 1}} f_1 \Argument{x_0 ^n} =
    \\
    & = x_0 ^{p_1 + 1} x_1 ^{p_2 - p_1} x_0 ^{p_3 - p_2} \dots x_{\Remainder{k
    + 1}} ^{p_k - p_{k - 1}} x_{1 - \Remainder{k + 1}} ^{\Argument{n - 1 - t_{k
    + 1} + p_{k + 1}} - p_{k}} x_{\Remainder{k + 1}} ^{n - 1 - \Argument{n - 1
    - t_{k + 1} + p_{k + 1}}},\\
    \s_2 \s_3 \Argument{x_0 ^n} & = f_0 f_1 \Argument{f_0 f_1} ^{p_1} f_1
    \Argument{f_0 f_1} ^{p_2} f_1 \dots \Argument{f_0 f_1} ^{p_k} f_1
    \Argument{f_0 f_1} ^{n - 1} = \\
    & = x_0 ^{p_1 + 1} x_1 ^{p_2 - p_1} x_0 ^{p_3 - p_2} \dots x_{1 -
    \Remainder{k}} ^{p_k - p_{k - 1}} x_{\Remainder{k}} ^{n - 1 - p_k},
\end{align*}
which contradicts Assumption~\eqref{eq:assumption}.

\par In the case $k > 0$ and $p_k \ge n - 1 - t_{k + 1} + p_{k + 1}$, let us
consider the element $\s_4 = \Argument{f_0 f_1} ^{n - 1 - t_{k + 1}} f_1
\Argument{f_0 f_1} ^{n - 1 - t_{k + 1} + p_{k + 1}}$. Then elements $\s_1 \s_4$
and $\s_2 \s_4$ are reduced to the following elements:
\begin{align*}
    \s_1 \s_4 & = f_0 f_1 \Argument{f_0 f_1} ^{p_1} f_1 \dots \Argument{f_0
    f_1} ^{p_k} f_1 \Argument{f_0 f_1} ^{p_{k + 1}} \cdot \Argument{f_0 f_1}
    ^{n - 1 - t_{k + 1}} f_1 \Argument{f_0 f_1} ^{n - 1 - t_{k + 1} + p_{k +
    1}} = \\
    & = f_0 f_1 \Argument{f_0 f_1} ^{p_1} f_1 \Argument{f_0 f_1} ^{p_2} f_1
    \dots \Argument{f_0 f_1} ^{p_{k - 1}} f_1 \Argument{f_0 f_1} ^{p_k};\\
    \s_2 \s_4 & = f_0 f_1 \Argument{f_0 f_1} ^{p_1} f_1 \dots \Argument{f_0
    f_1} ^{p_k} f_1 \Argument{f_0 f_1} ^{t_{k + 1}} \cdot \Argument{f_0 f_1}
    ^{n - 1 - t_{k + 1}} f_1 \Argument{f_0 f_1} ^{n - 1 - t_{k + 1} + p_{k +
    1}} = \\
    & = f_0 f_1 \Argument{f_0 f_1} ^{p_1} f_1 \Argument{f_0 f_1} ^{p_2} f_1
    \dots \Argument{f_0 f_1} ^{p_k} f_1 \Argument{f_0 f_1} ^{n - 1}.
\end{align*}
Similarly, for the input word $u = x_0 ^n$ we have
\begin{align*}
    \s_1 \s_4 \Argument{x_0 ^n} & = f_0 f_1 \Argument{f_0 f_1} ^{p_1} f_1
    \Argument{f_0 f_1} ^{p_2} f_1 \dots \Argument{f_0 f_1} ^{p_{k - 1}} f_1
    \Argument{f_0 f_1} ^{p_{k}} \Argument{x_0 ^n} = \\
    & = x_0 ^{p_1 + 1} x_1 ^{p_2 - p_1} x_0 ^{p_3 - p_2} \dots x_{1 -
    \Remainder{k - 1}} ^{p_{k - 1} - p_{k - 2}} x_{\Remainder{k -
    1}} ^{n - 1 - p_{k - 1}},\\
    \s_2 \s_4 \Argument{x_0 ^n} & = f_0 f_1 \Argument{f_0 f_1} ^{p_1} f_1
    \Argument{f_0 f_1} ^{p_2} f_1 \dots \Argument{f_0 f_1} ^{p_k} f_1
    \Argument{f_0 f_1} ^{n - 1} = \\
    & = x_0 ^{p_1 + 1} x_1 ^{p_2 - p_1} x_0 ^{p_3 - p_2} \dots x_{1 -
    \Remainder{k}} ^{p_k - p_{k - 1}} x_{\Remainder{k}} ^{n - 1 - p_k},
\end{align*}
which contradicts Assumption~\eqref{eq:assumption}. Hence, the
elements~\eqref{eq:elements_quotient} at $\epsilon _2 = \mu _2 = 0$
define the same transformation over $\Alphabet{2} ^n$ if and only if
$k = \ell$, $p_i = t_i$, $i = 1, 2, \dots, {k + 1}$.

\par Consider now the case $\epsilon _2 = 1$, $\mu _2 = 0$. Let us
assume that $\ell = 0$ or $t_\ell < t_{\ell + 1}$. In this case the elements
$\s_1 f_1$ and $\s_2 f_1$ are not reduced:
\begin{align*}
    \s_1 f_1 & = f_0 f_1 \Argument{f_0 f_1} ^{p_1} f_1 \Argument{f_0 f_1}
    ^{p_2} f_1 \dots \Argument{f_0 f_1} ^{p_k} f_1 \Argument{f_0 f_1} ^{p_{k +
    1} + 1},\\
    \s_2 f_1 & = f_0 f_1 \Argument{f_0 f_1} ^{t_1} f_1 \Argument{f_0 f_1}
    ^{t_2} f_1 \dots \Argument{f_0 f_1} ^{t_\ell} f_1 \Argument{f_0 f_1} ^{t_{\ell +
    1}} f_1.
\end{align*}
From Assumption~\eqref{eq:assumption} and the proof above the
requirements $k = \ell + 1$, $p_i = t_i$, $1 \le i \le k$, $p_{k + 1} + 1
= 0$ follow. The last requirement contradicts the condition $p_{k + 1}
\ge 0$ of~\eqref{eq:elements_quotient}. A similar reasoning can be
carried out for the elements $\s_1 f_0$ and $\s_2 f_0$, and we reach a
contradiction in the case $k = 0$ or $p_{k} < p_{k + 1}$.

\par Let us now consider the case $k,\ell > 0$, $t_\ell \ge t_{\ell + 1}$ and $p_{k} \ge
p_{k + 1}$. The elements $\s_1 f_1 \Argument{f_0 f_1} ^{t_{\ell + 1}}$ and $\s_2
f_1 \Argument{f_0 f_1} ^{t_{\ell + 1}}$ are reduced to the following normal forms:
\begin{align*}
    \s_1 f_1 \Argument{f_0 f_1} ^{t_{\ell + 1}} & = f_0 f_1 \Argument{f_0 f_1}
    ^{p_1} f_1 \Argument{f_0 f_1} ^{p_2} f_1 \dots \Argument{f_0 f_1} ^{p_k}
    f_1 \Argument{f_0 f_1} ^{\min\Pair{p_{k + 1} + 1 + t_{\ell + 1}}{n - 1}},\\
    \s_2 f_1 \Argument{f_0 f_1} ^{t_{\ell + 1}} & = f_0 f_1 \Argument{f_0 f_1}
    ^{t_1} f_1 \Argument{f_0 f_1} ^{t_2} f_1 \dots \Argument{f_0 f_1} ^{t_{\ell - 1}}
    f_1 \Argument{f_0 f_1} ^{t_{\ell}}.
\end{align*}
From Assumption~\eqref{eq:assumption} and the proof above the requirements
\begin{equation} \label{eq:requirements_quotient}
    k = \ell - 1, p_i = t_i, 1 \le i \le k, \min\Pair{p_{k + 1} + 1 + t_{\ell + 1}}{n
    - 1} = t_\ell
\end{equation}
follow. Similarly, from the equality \[
    \s_1 f_0 f_1 \Argument{f_0 f_1} ^{t_{\ell + 1}} \Argument{x_0 ^n} = \s_2 f_0
    f_1 \Argument{f_0 f_1} ^{t_{\ell + 1}} \Argument{x_0 ^n}
\] we get the requirement $k - 1 = \ell$, which contradicts the
requirements~\eqref{eq:requirements_quotient}.

\par The theorem is completely proved.
\end{proof}

\begin{proof}[Proof of Corollary~\ref{cor:Hausdorff_dimension}]
  Let us fix a number $n$, $n \ge 1$, and prove that the cardinality
  of the semigroup $\Quotient$ is \[ \Card{\Quotient} = 2 +
  \Argument{2n - 1} 2^n.
\]
Any element of form~\eqref{eq:normal_form_quotient} is defined by a set of $k$
parameters $\left\{ p_1, p_2, \dots, p_k \right\}$, and by $\epsilon_1$, $p_{k
+ 1}$, $\epsilon_2$. Parameter $\epsilon_1$ has two possible values, ``the
tail'' $\Argument{f_0 f_1} ^{p_{k + 1}} f_0 ^{\epsilon _2}$ has length varying
from $0$ to $\Argument{2n - 2}$, and the set $\left\{ p_1, p_2, \dots, p_k
\right\}$ is a $k$-element subset of $\left\{ 0, 1, \dots, n - 2 \right\}$,
where $k$ is some integer in $\left\{ 0, 1, \dots, n - 1 \right\}$. Therefore,
\[
    \Card{\Quotient} = \underbrace{2} _{1, f_0} + \underbrace{2}
    _{\epsilon_1} \cdot \underbrace{2 ^{n - 1}} _{p_1, \dots, p_k} \cdot
    \underbrace{\Argument{ 2n - 1 }} _{p_{k + 1} + \epsilon _2} = 2 +
    \Argument{ 2n - 1 } 2 ^n.
\]

\par As was shown  in Section~\ref{subsect:Hausdorff_dimension}, for all $n \ge
1$
 \[ \Card{\EndLevels{2}{n}} = 2^{2 \frac{2^n - 1}{2 - 1}} = 4 ^{2
   ^n - 1},
\] and the Hausdorff dimension of the semigroup $\Semigroup{\Intermediate{2}}$
is \[
    H\!\dim{\Semigroup{\Intermediate{2}}} = \liminf_{n \to \infty} \frac{\log
    \Argument{{2 + \Argument{2n - 1} 2 ^n}}}{\Argument{2 ^n - 1}\log 4} = 0,
\] which proves the corollary.
\end{proof}

\par Among groups, having Hausdorff dimension $0$ seems to be related to
being solvable; in~\cite{abert-v:dimension} it is shown that a solvable
group acting on a tree necessarily has dimension $0$. Since solvable
groups are all linear, this raises the question of whether
$\Semigroup{\Intermediate{2}}$ is linear. While this is probably not
the case, an even more fruitful analogy would be a linear embedding of
$\Semigroup{\Intermediate{2}}$ over a local ring, such that the
semigroups $\Quotient$ are the congruence quotients over that ring.

\section{Growth functions}
We derive, in this section, the growth series of the semigroup
$\Semigroup{\Intermediate{2}}$, as well as the asymptotics of the growth
functions $\GrowthSemigroup{\Intermediate{2}}$ and $\Growth{\Intermediate{2}}$.

\subsection{Growth series}

\begin{lemma}\label{lemma:q}
  Let $\OddPartitions{n}$ be the number of partitions of $n \in \Natural$ in
  distinct, odd parts, and form $\Psi \Argument{X} = \sum \OddPartitions{n} X
  ^n$. Then
  \[\Psi \Argument{X} = \sum \limits _{m = 0} ^\infty \frac{X ^{m ^2}}{\left( 1
  - X ^2 \right) \cdots \left(1 - X ^{2m} \right)} = \left( 1 + X \right)
  \left( 1 + X ^3 \right) \left( 1 + X ^5 \right) \cdots\]
\end{lemma}

\begin{proof}
  Let $\lambda = \left( \lambda_1, \lambda_2, \dots, \lambda_m
  \right)$ be such a partition of $n$. Then $\lambda_i \ge 2i - 1$ for
  all $i = 1, 2, \dots, m$, and \[\left( \lambda_1 - 1, \lambda_2 - 3,
    \dots, \lambda_m - \left( 2m - 1 \right) \right)\] is a partition
  of $n - m ^2$ in at most $m$ even parts. By ``flipping'', this is
  the same as a partition of $n - m ^2$ into even parts that are at
  most $2m$, whence the first equality.

  \par The second equality is standard: an integer partition $\left( \lambda_1,
  \dots, \lambda_m \right)$ in distinct odd parts corresponds to a monomial $X
  ^{\lambda_1} \dots X ^{\lambda_m}$.
\end{proof}

\par It follows from Proposition~\ref{prop:normal_form} that the word growth
series of $\Semigroup{\Intermediate{2}}$ is
\[\Delta _{\Semigroup{\Intermediate{2}}} \Argument{X} = \sum \limits_{n \ge 0} {\WordGrowthSemigroup{\Intermediate{2}}
\n X ^n} = \underbrace{\Argument{1 + X}} _{1, f_0} + \underbrace{ \underbrace{
\Argument{1 + X}} _{f_0 ^{\epsilon_1}} \underbrace{X} _{f_1} \Psi \Argument{X}
\underbrace{\frac{1}{1 - X ^2}} _{\Argument{f_0 f_1} ^{p_{k + 1}}}
\underbrace{\Argument{1 + X}} _{f_1 ^{\epsilon_2}}} _{\text{form~\eqref{eq:normal_form}}}.
\]

\par Indeed all subwords $ \Argument{f_0 f_1} ^{p_1} f_1 \Argument{f_0 f_1}
^{p_2} \dots \Argument{f_0 f_1} ^{p_k} f_1$ of the second form~\eqref{eq:normal_form} correspond uniquely to an integer partition
$\left( 2p_1 + 1, 2p_2 + 1, \dots, 2p_{k} + 1 \right)$ in distinct
odd parts; and $\Argument{2p_1 + 1} + \Argument{2p_2 + 1} + \dots +
\Argument{2p_k + 1}$ is the length of this subword. We obtain:
\begin{equation} \label{eq:word_growth_series}
\begin{aligned}
  \Delta _{\Semigroup{\Intermediate{2}}} \n[X] & = 1 + X + \frac{X + X^2}{1 -
  X} \Psi \n[X] = \Argument{1 + X} \left( 1 + \frac{X}{1 - X} \prod \limits _{n
  \ge 0} \Argument{1 + X ^{2n + 1}} \right) = \\
  & = \Argument{1 + X} \left( 1 + \frac{X}{1 - X} \left( 1 + \frac{X}{1 - X^2}
  \left( 1 + \frac{X^3}{1 - X^4} \left( 1 + \dots \right) \right) \right)
  \right), \\
\end{aligned}
\end{equation}
that proves the first part of Theorem~\ref{th:generating_functions}.

\par As mentioned in Remark~\ref{rem:even_reducing}, the set of
elements which can be presented as a product of $n$ generators,
includes the sets of elements of length $n,n - 2,\dots$.  Therefore
\[
    \Growth{\Intermediate{2}} \n = \sum \limits _{i = 0} ^{\Divider{n}}
    {\WordGrowthSemigroup{\Intermediate{2}} \Argument{2i + \Remainder{n}}}
\] whence \[
    \Gamma _{\Intermediate{2}} \Argument{X} = \frac{1}{1 - X^2} \Delta
    _{\Semigroup{\Intermediate{2}}} \Argument{X} = \frac{1}{1 - X} \Argument{1
    + \frac{X}{1 - X} \prod \limits _{n \ge 0} \Argument{1 + X ^{2n + 1}}}.
\] As $\GrowthSemigroup{\Intermediate{2}} \n = \sum \limits _{i =
0} ^n {\WordGrowthSemigroup{\Intermediate{2}} \Argument{i}}$, one has \[
    \Gamma _{\Semigroup{\Intermediate{2}}} \Argument{X} = \frac{1}{1 - X}
    \Delta _{\Semigroup{\Intermediate{2}}} \Argument{X} = \frac{1 + X}{1 - X}
    \Argument{1 + \frac{X}{1 - X} \prod \limits _{n \ge 0} \Argument{1 + X ^{2n
    + 1}}}.
\] Last two equalities complete the proof of Theorem~\ref{th:generating_functions}.

\subsection{Asymptotics}

\par We quote the following result by Richmond~\cite{richmond:andrews}:
\begin{theorem}
  If $\gcd(a_1,\cdots,a_s,M)=1$, then the number of partitions of $n$
  into distinct parts all congruent to some $a_i \mod M$ has the
  asymptotic value
\[
    2 ^{\left( \frac{s - 3}{2} + \frac{1}{M}\left( \sum a_i \right) \right)} 3
    ^{-1/4} n ^{-3/4} \exp \left(\pi \sqrt{\frac{sn}{3M}} \right) \left( 1 +
    \mathcal{O} \left( n ^{-1/2 + \delta} \right) \right)
\] for any $\delta > 0$.
\end{theorem}

\par In particular, for $\OddPartitions{n}$, we take $M = 2$, $s = 1$ and
$a_1 = 1$ to obtain the asymptotics
\[
    \OddPartitions{n} \sim 2^{-1/2} 3^{-1/4} n^{-3/4} \exp \Argument{\pi
    \sqrt{\frac{n}{6}}},
\] where $f \n \sim g \n$ means $\lim {f \n} / {g \n} = 1$.

The following result appears as Lemma~3.4
in~\cite{lavrik:interm}. Its proof follows from the
Euler-MacLaurin summation formula:
\begin{lemma}[\cite{lavrik:interm}]\label{lemma:sum}
  Let $f$ be a series with $f \n \sim n ^\alpha \exp
  \Argument{\beta\sqrt n}$, and define $g \n = \sum \limits _{i = 1}
  ^n f \n[i]$. Then \[g \n \sim \frac{2}{\beta} n ^{\alpha + 1/2} \exp
  \Argument{\beta \sqrt n}.\]
\end{lemma}


\par Let us return to the first expression in~\eqref{eq:word_growth_series}.
The term $\Argument{X + X^2} / \Argument{1 - X}$ expands to $X + 2X ^2 + 2X ^3
+ \dots$. We deduce:
\begin{subequations}
\begin{equation} \label{eq:word_growth_partitions}
    \WordGrowthSemigroup{\Intermediate{2}} \n = \OddPartitions{n - 1} + 2 \sum
    \limits _{i = 0} ^{n - 2} \OddPartitions{i}
\end{equation}
for $n \ge 2$. Moreover,
\begin{align}
    \label{eq:growth_automaton_partitions}
    \Growth{\Intermediate{2}} \n & = 1 + \sum \limits _{i = 0} ^{n - 1}
    \Argument{n - i} \OddPartitions{i},\\
    \label{eq:growth_semigroup_partitions}
    \GrowthSemigroup{\Intermediate{2}} \n & = 2 + \sum \limits _{i = 0} ^{n -
    1} \Argument{2n - 2i - 1} \OddPartitions{i}.
\end{align}
\end{subequations}


\begin{proof}[Proof of Theorem~\ref{th:estimates}]
It follows from Lemma~\ref{lemma:sum} and~\eqref{eq:word_growth_partitions}
that \[
    \WordGrowthSemigroup{\Intermediate{2}} \n \sim 2 \sum \limits _{i = 0} ^{n}
    {\OddPartitions{i}} \sim \frac{4 \sqrt{6}}{\pi}
    \sqrt{n} \cdot \OddPartitions{n} \sim \frac{2 ^2 3 ^{1/4}}{\pi} n^{-1/4}
    \exp{\Argument{\pi \sqrt{\frac{n}{6}}}}.
\]

\par Once more, applying Lemma~\ref{lemma:sum} to the equation at line above,
we have the sharp estimate \[
    \GrowthSemigroup{\Intermediate{2}} \n = \sum \limits _{i = 0} ^{n}
    \WordGrowthSemigroup{\Intermediate{2}} \n[i] \sim \frac{48}{\pi ^2} n \cdot
    \OddPartitions{n} \sim \frac{2 ^{7/2} 3 ^{3/4}}{\pi ^2} n^{1/4}
    \exp{\Argument{\pi \sqrt{\frac{n}{6}}}},
\]
with the ratios of left- to right-hand side tending to $1$ as
$n \to \infty$.

\par Similarly, the growth function of the automaton
$\Intermediate{2}$ admits the sharp estimate \[
    \Growth{\Intermediate{2}} \n \sim \frac{24}{\pi ^2} n \cdot
    \OddPartitions{n} \sim \frac{2 ^{5/2} 3 ^{3/4}}{\pi ^2} n^{1/4}
    \exp{\Argument{\pi \sqrt{\frac{n}{6}}}},
\] that completes the proof of Theorem~\ref{th:estimates}.
\end{proof}

\begin{proof}[Proof of Corollary~\ref{cor:growth_orders}]
  From Theorem~\ref{th:estimates} it follows that \[
  \GrowthOrder{\GrowthSemigroup{\Intermediate{2}}} = \GrowthOrder{\exp
    \n[\sqrt n]},
\] and by Proposition~\ref{prop:monoid_growth_order} the same asymptotics hold for
$\GrowthOrder{\Growth{\Intermediate{2}}}$.
\end{proof}

\end{document}